\documentclass[11pt,a4paper,leqno]{amsart}
\usepackage{amssymb}
\usepackage[noPostScript=dvips]{diagrams}
\newif\ifaddpics\addpicstrue
\swapnumbers
\theoremstyle{plain}
\newtheorem*{result}{Theorem}
\newtheorem{thm}{Theorem}[section]
\newtheorem{prop}[thm]{Proposition}

\theoremstyle{definition}
\newtheorem{defn}[thm]{Definition}
\newtheorem{conv}[thm]{Conventions}
\theoremstyle{remark}
\newtheorem{rem}[thm]{Remark}
\newcommand{\proofof}[1]{\end{#1}\begin{proof}}
\newcommand{\emphdef}{\textit}
\newcommand{\ip}[1]{\langle#1\rangle}

\newcommand{\R}{{\mathbb R}}
\newcommand{\C}{{\mathbb C}}

\newcommand{\tens}{\mathbin{\otimes}}
\newcommand{\skwend}{\mathinner{\vartriangle}}

\newcommand{\setdif}{\smallsetminus}
\newcommand{\dual}{^{*\!}}
\newcommand{\Cinf}{\mathrm{C}^\infty}
\newcommand{\Mn}{M}

\newcommand{\SO}{\mathrm{SO}}

\newcommand{\GL}{\mathrm{GL}}
\newcommand{\SL}{\mathrm{SL}}
\newcommand{\Un}{\mathrm{U}}
\newcommand{\SU}{\mathrm{SU}}

\DeclareMathAlphabet{\mathrmsl}{OT1}{cmr}{m}{sl}

\newcommand{\rssymb}[2]{\newcommand{#1}{\mathrmsl{#2}} }
\newcommand{\oper}[3][n]{\newcommand{#2}{\mathop{\mathrm{#3}}%
\ifx n#1\nolimits\else\limits\fi} }
\newcommand{\rsoper}[3][n]{\newcommand{#2}{\mathop{\mathrmsl{#3}}%
\ifx n#1\nolimits\else\limits\fi} }
\oper\End{End}
\oper\Aut{Aut}
\oper\Alt{Alt}
\oper\Sym{Sym}
\rsoper\divg{div}
\rsoper{\sym}{sym}
\rsoper{\alt}{alt}
\rsoper\trace{tr}
\rssymb\iden{id}
\newcommand{\acknowledge}{\subsection*{Acknowledgements.}}
\newcommand{\thismonth}{\ifcase\month\or
  January\or February\or March\or April\or May\or June\or
  July\or August\or September\or October\or November\or December\fi
  \space\number\year}
\newcommand{\tcite}[1]{~\textup{\cite{#1}}}
\newcommand{\bauth}[1]{\mbox{#1},}
\newcommand{\bart}[1]{\textit{#1},}
\newcommand{\bjourn}[3]{#1 \textbf{#2} (#3)}
\newcommand{\bbook}[1]{\textsl{#1},}

\newcommand{\bpp}[2]{pp.~#1--#2.}
\numberwithin{equation}{section}

\newcommand{\conf}{\mathsf{c}}
\newcommand{\cip}{\ip}
\newcommand{\mult}{^{\scriptscriptstyle\times}}
\rsoper{\Ric}{Ric}
\rsoper{\ric}{ric}
\rssymb{\scal}{scal}

\let\geq\geqslant

\let\connect\#
\def\s/{^{\vphantom{,}}_}
\def\S/{_{\vphantom{T}}^}
\newcommand{\D}{D^{\vphantom{x}}}
\newcommand{\sd}{{\scriptscriptstyle +}}
\newcommand{\asd}{{\scriptscriptstyle -}}
\newcommand{\sdasd}{{\scriptscriptstyle\pm}}
\newcommand{\CP}[1]{\C P^{#1}}

\newcommand{\cD}{\mathcal D}
\newcommand{\cC}{\mathcal C}
\newcommand{\cL}{\mathcal L}
\newcommand{\cH}{\mathcal H}
\newcommand{\cO}{\mathcal O}
\newcommand{\cS}{{\mathcal S}}
\newcommand{\gmw}{\mathrmsl w}

\newcommand{\dv}{\tau}
\newcommand{\dvK}{\dv\s/0}
\newcommand{\tw}{\kappa}
\newcommand{\twK}{\tw\s/0}
\newcommand{\sh}{\sigma}
\newcommand{\Sh}{\Sigma}
\newcommand{\dr}{\rho}
\newcommand{\ortn}{\underline{\mathrmsl{or}}}
\ifaddpics\input epsf
\newcommand{\getpic}[1]{\epsfxsize=.75truein\lower.3truein\hbox
{\epsfbox{#1}}}
\newcommand{\getbigpic}[1]{\epsfxsize=1.8truein\epsfbox{#1}}\fi
\begin{document}
\title[Selfdual spaces and Einstein-Weyl geometry]{Selfdual spaces
with complex structures,\\
Einstein-Weyl geometry and geodesics}
\author{David M. J. Calderbank}
\address{Department of Mathematics and Statistics\\
University of Edinburgh\\ King's Buildings, Mayfield Road\\
Edinburgh EH9 3JZ\\ Scotland.}
\email{davidmjc@maths.ed.ac.uk}
\author{Henrik Pedersen}
\address{Department of Mathematics and Computer Science\\
Odense University\\ Campusvej 55\\ DK-5230 Odense M\\ Denmark.}
\email{henrik@imada.ou.dk}
\date{\thismonth}
\begin{abstract}
  We study the Jones and Tod correspondence between selfdual conformal
  $4$-manifolds with a conformal vector field and abelian monopoles on
  Einstein-Weyl $3$-manifolds, and prove that invariant complex structures
  correspond to shear-free geodesic congruences. Such congruences exist in
  abundance and so provide a tool for constructing interesting selfdual
  geometries with symmetry, unifying the theories of scalar-flat K\"ahler
  metrics and hypercomplex structures with symmetry. We also show that
  in the presence of such a congruence, the Einstein-Weyl equation is
  equivalent to a pair of coupled monopole equations, and we solve these
  equations in a special case. The new Einstein-Weyl spaces, which
  we call Einstein-Weyl ``with a geodesic symmetry'', give rise to
  hypercomplex structures with two commuting triholomorphic vector fields.
\end{abstract}
\maketitle

\section{Introduction}

Selfdual conformal $4$-manifolds play a central role in low dimensional
differential geometry. The selfduality equation is integrable, in the sense
that there is a twistor construction for solutions, and so one can hope to
find many explicit examples~\cite{Besse,MW}. One approach is to look for
examples with symmetry. Since the selfduality equation is the complete
integrability condition for the local existence of orthogonal (and
antiselfdual) complex structures, it is also natural to look for solutions
equipped with such complex structures. Our aim herein is to study the geometry
of this situation in detail and present a framework unifying the theories of
hypercomplex structures and scalar-flat K\"ahler metrics with
symmetry~\cite{CTV,GT,LeBrun1}. Within this framework, there are explicit
examples of hyperK\"ahler, selfdual Einstein, hypercomplex and scalar-flat
K\"ahler metrics parameterised by arbitrary functions.

The key tool in our study is the Jones and Tod construction~\cite{JT},
which shows that the reduction of the selfduality equation by
a conformal vector field is given by the Einstein-Weyl equation together
with the linear equation for an abelian monopole. This correspondence
between a selfdual space $M$ with symmetry and an Einstein-Weyl space $B$
with a monopole is remarkable for three reasons:
\begin{enumerate}
\item It provides a geometric interpretation of the symmetry reduced equation
for an arbitrary conformal vector field.
\item It is a constructive method for building selfdual spaces out of
solutions to a linear equation on an Einstein-Weyl space.
\item It can be used in the other direction to construct Einstein-Weyl spaces
from selfdual spaces with symmetry.
\end{enumerate}

We add to this correspondence by proving that invariant antiselfdual complex
structures on $M$ correspond to shear-free geodesic congruences on $B$, i.e.,
foliations of $B$ by oriented geodesics, such that the transverse conformal
structure is invariant along the leaves.  This generalises Tod's
observation~\cite{Tod2} that the Einstein-Weyl spaces arising from scalar-flat
K\"ahler metrics with Killing fields~\cite{LeBrun1} admit a shear-free geodesic
congruence which is also twist-free (i.e., surface-orthogonal).

In order to explain how the scalar-flat K\"ahler story and the analogous story
for hypercomplex structures~\cite{CTV,GT} fit into our more general narrative,
we begin, in section~\ref{CSKW}, by reviewing, in a novel way, the
construction of a canonical ``K\"ahler-Weyl connection'' on any conformal
Hermitian surface~\cite{Gauduchon1,Vaisman1}. We give a representation
theoretic proof of the formula for the antiselfdual Weyl tensor on such a
surface~\cite{AG} and discuss its geometric and twistorial interpretation when
the antiselfdual Weyl tensor vanishes. We use twistor theory throughout the
paper to explain and motivate the geometric constructions, although we
find it easier to make these constructions more general, explicit and precise
by direct geometric arguments.

Having described the four dimensional context, we lay the three dimensional
foundations for our study in section~\ref{SGCEW}. We begin with some
elementary facts about congruences, and then go on to show that the
Einstein-Weyl equation is the complete integrability condition for the
existence of shear-free geodesic congruences in a three dimensional Weyl
space. As in section~\ref{CSKW}, we discuss the twistorial interpretation,
this time in terms of the associated ``minitwistor space''~\cite{Hitchin3},
and explain the minitwistor version of the Kerr theorem, which has only been
discussed informally in the existing literature (and usually only in the flat
case). We also show that at any point where the Einstein-Weyl condition does
not hold, there are at most two possible directions for a shear-free geodesic
congruence. The main result of our work in this section, however, is a
reformulation of the Einstein-Weyl equation in the presence of a shear-free
geodesic congruence. More precisely, we show in Theorem~\ref{EWeqnthm} that
the Einstein-Weyl equation is equivalent to the fact that the divergence and
twist of this congruence are both monopoles of a special kind. These
``special'' monopoles play a crucial role in the sequel.

We end section~\ref{SGCEW} by giving examples. We first explain how the
Einstein-Weyl spaces arising as quotients of scalar-flat K\"ahler metrics and
hypercomplex structures fit into our theory: they are the cases of vanishing
twist and divergence respectively. In these cases it is known that the
remaining nonzero special monopole (i.e., the divergence and twist
respectively) may be used to construct a hyperK\"ahler
metric~\cite{BF,CTV,GT}, motivating some of our later results. We also give
some new examples: indeed, in Theorem~\ref{gs}, we classify explicitly the
Einstein-Weyl spaces admitting a geodesic congruence generated by a conformal
vector field preserving the Weyl connection. We call such spaces Einstein-Weyl
\emphdef{with a geodesic symmetry}. They are parameterised by an arbitrary
holomorphic function of one variable.

The following section contains the central results of this paper, in which the
four and three dimensional geometries are related.  We begin by giving a new
differential geometric proof of the Jones and Tod correspondence~\cite{JT}
between oriented conformal structures and Weyl structures, which reduces the
selfduality condition to the Einstein-Weyl condition
(see~\ref{JonesnTod}). Although other direct proofs can be found in the
literature~\cite{GT,Joyce,LeBrun1}, they either only cover special cases, or
are not sufficiently explicit for our purposes. Our next result,
Theorem~\ref{csthm}, like the Jones and Tod construction, is motivated by
twistor theory. Loosely stated, it is as follows.
\begin{result} Suppose $M$ is an oriented conformal $4$-manifold with a
conformal vector field, and $B$ is the corresponding Weyl space. Then
invariant antiselfdual complex structures on $M$ correspond to shear-free
geodesic congruences on $B$.
\end{result}

In fact we show explicitly how the K\"ahler-Weyl connection may be constructed
from the divergence and twist of the congruence. This allows us to
characterise the hypercomplex and scalar-flat K\"ahler cases of our
correspondence, reobtaining the basic constructions
of~\cite{BF,CTV,GT,LeBrun1}, as well as treating quotients of hypercomplex,
scalar-flat K\"ahler and hyperK\"ahler manifolds by more general holomorphic
conformal vector fields. As a consequence, we show in Theorem~\ref{quotthm}
that every Einstein-Weyl space is locally the quotient of some scalar-flat
K\"ahler metric and also of some hypercomplex structure, and that it is a
local quotient of a hyperK\"ahler metric (by a holomorphic conformal vector
field) if and only if it admits a shear-free geodesic congruence with linearly
dependent divergence and twist.

We clarify the scope of these results in section~\ref{SDE} where we show that
our constructions can be applied to all selfdual Einstein metrics with a
conformal vector field. Here, we make use of the fact that a selfdual Einstein
metric with a Killing field is conformal to a scalar-flat K\"ahler
metric~\cite{Tod3}.

The last four sections are concerned exclusively with examples.  In
section~\ref{EWR4} we show how our methods provide some insight into the
construction of Einstein-Weyl structures from $\R^4$~\cite{PT1}. As a
consequence, we observe that there is a one parameter family of Einstein-Weyl
structures on $S^3$ admitting shear-free twist-free geodesic congruences.
This family is complementary to the more familiar Berger spheres, which admit
shear-free divergence-free geodesic congruences~\cite{CTV,GT}.

In section~\ref{HKt}, we generalise this by replacing $\R^4$ with a
Gibbons-Hawking hyperK\"ahler metric~\cite{GH} constructed from a harmonic
function on $\R^3$. If the corresponding monopole is invariant under a
homothetic vector field on $\R^3$, then the hyperK\"ahler metric has an extra
symmetry, and hence another quotient Einstein-Weyl space. We first treat the
case of axial symmetry, introduced by Ward~\cite{Ward}, and then turn to more
general symmetries. The Gibbons-Hawking metrics constructed from monopoles
invariant under a general Killing field give new implicit solutions of the
Toda field equation. On the other hand, from the monopoles invariant under
dilation, we reobtain the Einstein-Weyl spaces with geodesic symmetry.

In section~\ref{CMCC} we look at the constant curvature metrics on $\cH^3$,
$\R^3$ and $S^3$ from the point of view of congruences and use this prism to
explain properties of the selfdual Einstein metrics fibering over them. Then
in the final section, we consider once more the Einstein-Weyl spaces
constructed from harmonic functions on $\R^3$, and use them to construct torus
symmetric selfdual conformal structures. These include those of
Joyce~\cite{Joyce}, some of which live on $k\CP2$, and also an explicit family
of hypercomplex structures depending on two holomorphic functions of one
variable.

This paper is primarily concerned with the richness of the local geometry of
selfdual spaces with symmetry, and we have not studied completeness or
compactness questions in any detail. Indeed, the local nature of the Jones and
Tod construction makes it technically difficult to tackle such issues from
this point of view, and doing so would have added considerably to the length
of this paper. Nevertheless, there remain interesting problems which we hope
to address in the future.

\acknowledge Thanks to Paul Gauduchon, Michael Singer and Paul Tod for helpful
discussions. The diagrams were produced using Xfig, Mathematica and Paul
Taylor's commutative diagrams package.

\section{Conformal structures and K\"ahler-Weyl geometry}\label{CSKW}

Associated to an orthogonal complex structure $J$ on a conformal manifold is a
distinguished torsion-free connection $D$. The conformal structure is
preserved by this connection and, in four dimensions, so is $J$. Such a
connection is called a \emphdef{K\"ahler-Weyl connection}~\cite{CP}: if it is
the Levi-Civita connection of a compatible Riemannian metric, then this metric
is K\"ahler. In this section, we review this construction, which goes back to
Lee and Vaisman (see~\cite{Gauduchon1,Lee,Vaisman1}).

It is convenient in conformal geometry to make use of the \emphdef{density
bundles} $L^w$ (for $w\in\R$). On an $n$-manifold $\Mn$, $L^w$ is the oriented
real line bundle associated to the frame bundle by the representation
$A\mapsto|\det A|^{w/n}$ of $\GL(n)$. The fibre $L^w_x$ may be constructed
canonically as the space of maps $\rho\colon(\Lambda^nT_x\Mn)\setdif0\to\R$
such that $\rho(\lambda\omega)=|\lambda|^{-w/n}\rho(\omega)$ for all
$\lambda\in\R\mult$ and $\omega\in(\Lambda^nT_x\Mn)\setdif0$.

A \emphdef{conformal structure} $\conf$ on $\Mn$ is a positive definite
symmetric bilinear form on $T\Mn$ with values in $L^2$, or equivalently a
metric on the bundle $L^{-1}T\Mn$. (When tensoring with a density line
bundle, we generally omit the tensor product sign.)

The line bundles $L^w$ are trivialisable and a nonvanishing (usually
positive) section of $L^1$ (or $L^w$ for $w\neq0$) will be called a
\emphdef{length scale} or \emphdef{gauge} (of weight $w$).  We also say that
tensors in $L^w\tens(T\Mn)^j\tens(T\dual\Mn)^k$ have \emphdef{weight}
$w+j-k$.  If $\mu$ is a positive section of $L^1$, then $\mu^{-2}\conf$ is a
Riemannian metric on $\Mn$, which will be called~\emphdef{compatible}.  A
conformal structure may equally be defined by the associated ``conformal
class'' of compatible Riemannian metrics.

A \emphdef{Weyl derivative} is a covariant derivative $D$ on $L^1$. It
induces covariant derivatives on $L^w$ for all $w$.  The curvature of $D$ is
a real $2$-form $F^D$ which will be called the \emphdef{Faraday curvature} or
\emphdef{Faraday $2$-form}.  If $F^D=0$ then $D$ is said to be
\emphdef{closed}. It follows that there are local length scales $\mu$ with
$D\mu=0$. If such a length scale exists globally then $D$ is said to be
\emphdef{exact}.  Conversely, a length scale $\mu$ induces an exact Weyl
derivative $D^\mu$ such that $D^\mu\mu=0$. Consequently, we sometimes refer
to an exact Weyl derivative as a \emphdef{gauge}. The space of Weyl
derivatives on $\Mn$ is an affine space modelled on the space of $1$-forms.

Any connection on $T\Mn$ induces a Weyl derivative on $L^1$. Conversely,
on a conformal manifold, the Koszul formula shows that any Weyl derivative
determines uniquely a torsion-free connection $D$ on $T\Mn$ with $Dc=0$
(see~\cite{CP}). Such connections are called \emphdef{Weyl
connections}. Linearising the Koszul formula with respect to $D$ shows
that $(D+\gamma)_XY=D_XY+\gamma(X)Y+\gamma(Y)X-\cip{X,Y}\gamma$,
where $\cip{.\,,.}$ denotes the conformal structure, and $X,Y$ are
vector fields. Notice that here, and elsewhere, we make free use of the sharp
isomorphism $\sharp\colon T\dual\Mn\to L^{-2}T\Mn$. We sometimes
write $\gamma\skwend X (Y)=\iota\s/Y(\gamma\wedge X)$ for the last
two terms.

\begin{defn} A \emphdef{K\"ahler-Weyl} structure on a conformal manifold
$\Mn$ is given by a Weyl derivative $D$ and an orthogonal
complex structure $J$ such that $DJ=0$.
\end{defn}

Suppose now that $\Mn$ is a conformal $n$-manifold ($n=2m>2$) and that $J$ is
an orthogonal complex structure. Then $\Omega_J:=\cip{J.,.}$ is a section
of $L^2\Lambda^2T\dual\Mn$, called the \emphdef{conformal K\"ahler form}.
It is a nondegenerate weightless $2$-form. [In general, we identify
bilinear forms and endomorphism by $\Phi(X,Y)=\cip{\Phi(X),Y}$.]

\begin{prop}\textup{(cf.~\cite{Lee})} Suppose that $\Omega$ is a nondegenerate
weightless $2$-form. Then there is a unique Weyl derivative $D$ such that
$d^D\Omega$ is trace-free with respect to $\Omega$, in the sense that $\sum
d^D\Omega(e_i,e_i',.)=0$, where $e_i,e_i'$ are frames for $L^{-1}T\Mn$ with
$\Omega(e_i,e_j')=\delta_{ij}$.  \proofof{prop} Pick any Weyl derivative $D^0$
and set $D=D^0+\gamma$ for some $1$-form $\gamma$. Then
$d^D\Omega=d^{D^0}\Omega+2\gamma\wedge\Omega$ and so the traces differ by
\begin{align*}
2(\gamma\wedge\Omega)(e_i,e_i',.)
&=2\gamma(e_i)\Omega(e_i',.)+2\gamma(e_i')\Omega(.,e_i)+
2\gamma\,\Omega(e_i,e_i')\\
&=2(n-2)\gamma.
\end{align*}
Since $n>2$ it follows that there is a unique $\gamma$ such that
$d^D\Omega$ is trace-free.
\end{proof}

\begin{prop}\label{KWprop} Suppose that $J$ is an orthogonal complex
structure on a conformal manifold $\Mn$ and that $d^D\Omega_J=0$. Then $D$
defines a K\"ahler-Weyl structure on $\Mn$, i.e., $DJ=0$.

\proofof{prop} For any vector field $X$, $D_XJ$ anticommutes
with $J$ (since $J^2=-\iden$) and is skew (since $J$ is skew, and
$D$ is conformal). Hence $\cip{(D_{JX}J-JD_XJ)Y,Z}$, which is
symmetric in $X,Y$ because $J$ is integrable and $D$ is torsion-free,
is also skew in $Y,Z$. It must therefore vanish for all $X,Y,Z$.
If we now impose $d^D\Omega_J=0$ we obtain:
\begin{align*}
0&=d^D\Omega(X,Y,Z)-d^D\Omega(X,JY,JZ)\\
&=\cip{(D_XJ)Y,Z}+\cip{(D_YJ)Z,X}+\cip{(D_ZJ)X,Y}\\
&\quad-\cip{(D_XJ)JY,JZ}-\cip{(JD_YJ)JZ,X}-\cip{(JD_ZJ)X,JY}\\
&=2\cip{(D_XJ)Y,Z}.
\end{align*}
Hence $DJ=0$.
\end{proof}

Now if $n=4$ and $D$ is the unique Weyl derivative such that
$d^D\Omega_J$ is trace-free, then in fact $d^D\Omega_J=0$ since
wedge product with $\Omega_J$ is an isomorphism from $T\dual\Mn$ to
$L^2\Lambda^3T\dual\Mn$. Hence, by Proposition~\ref{KWprop}, $DJ=0$.
To summarise:
\begin{thm}\label{KWthm}\tcite{Vaisman1} Any Hermitian conformal
structure on any complex surface $\Mn$ induces a unique K\"ahler-Weyl
structure on $\Mn$. The Weyl derivative is exact iff the conformal Hermitian
structure admits a compatible K\"ahler metric.
\end{thm}

On an oriented conformal $4$-manifold, orthogonal complex structures are
either selfdual or antiselfdual, in the sense that the conformal K\"ahler form
is either a selfdual or an antiselfdual weightless $2$-form. In this paper we
shall be concerned primarily with antiselfdual complex structures on selfdual
conformal manifolds, i.e., conformal manifolds $\Mn$ with $W^\asd=0$, where
$W^\asd$ is the antiselfdual part of the Weyl tensor. In this case, as is well
known (see~\cite{Besse}), there is a complex $3$-manifold $Z$ fibering over
$\Mn$, called the \emphdef{twistor space} of $\Mn$. The fibre $Z_x$ given by
the $2$-sphere of orthogonal antiselfdual complex structures on $T_x\Mn$, and
the antipodal map $J\mapsto -J$ is a real structure on $Z$. The fibres are
called the (real) \emphdef{twistor lines} of $Z$ and are holomorphic rational
curves in $Z$. The canonical bundle $K_Z$ of $Z$ is easily seen to be of
degree $-4$ on each twistor line. As shown in~\cite{Gauduchon2,PS1}, any Weyl
derivative on $\Mn$ whose Faraday $2$-form is selfdual induces a holomorphic
structure on $L^1_\C$, the pullback of $L^1\tens\C$, and (up to reality
conditions) this process is invertible; this is the Ward correspondence for
line bundles, or the Penrose correspondence for selfdual Maxwell fields.

The K\"ahler-Weyl connection arising in Theorem~\ref{KWthm} can be given a
twistor space interpretation. Any antiselfdual complex structure $J$ defines
divisors $\cD,\overline\cD$ in $Z$, namely the sections of $Z$ given by
$J,-J$. Since the divisor $\cD+\overline\cD$ intersects each twistor line
twice, the holomorphic line bundle $[\cD+\overline\cD]K_Z^{1/2}$ is trivial on
each twistor line: more precisely, by viewing $J_x$ as a constant vector field
on $L^2\Lambda^2T\dual_x M$, its orthogonal projection canonically defines a
vertical vector field on $Z$ holomorphic on each fibre and vanishing along
$\cD+\overline\cD$. Therefore $[\cD+\overline\cD]$ is a holomorphic structure
on the vertical tangent bundle of $Z$. In fact the vertical bundle of $Z$ is
$L^{-1}_\C K_Z^{-1/2}$ and so $J$ determines a holomorphic structure on
$L^{-1}_\C$, which, since $[\cD+\overline\cD]$ is real, gives a Weyl
derivative on $\Mn$ with selfdual Faraday curvature~\cite{Gauduchon4}.

Similarly, by projecting each twistor line stereographically onto the
orthogonal complement of $J$ in $L^2\Lambda^2_\asd T\dual\Mn$, which we denote
$L^2K_J$, we see that the pullback of $L^2K_J$ to $Z$ has a section $s$
meromorphic on each fibre with a zero at $J$ and a pole at $-J$. Therefore the
divisor $\cD-\overline\cD$ defines a holomorphic structure on this pullback
bundle and hence a covariant derivative with (imaginary) selfdual curvature on
$L^2K_J$. This curvature may be identified with the \emphdef{Ricci form},
since if it vanishes, $[\cD-\overline\cD]$ is trivial, and so $s$, viewed as a
meromorphic function on $Z$, defines a fibration of $Z$ over $\CP1$;
that is, $M$ is hypercomplex.

The selfduality of the Faraday and Ricci forms may be deduced directly from
the selfduality of the Weyl tensor. To see this, we need a few basic facts
from Weyl and K\"ahler-Weyl geometry.

First of all, let $D$ be a Weyl derivative on a conformal $n$-manifold
and let $R^{D,w}$ denote the curvature of $D$ on $L^{w-1}T\Mn$.
Then it is well known that:
\begin{equation}\label{curv} R^{D,w}_{X,Y}=W\s/{X,Y}+wF^D(X,Y)\iden
-r^D(X)\skwend Y+r^D(Y)\skwend X.
\end{equation}
Here $W$ is the Weyl tensor and $r^D$ is the normalised Ricci tensor, which
decomposes under the orthogonal group as
$r^D=r^D_0+\frac1{2n(n-1)}\scal^D\iden-\frac12F^D$, where $r^D_0$ is
symmetric and trace-free, and the trace part defines the scalar curvature of
$D$.

\begin{prop} On a K\"ahler-Weyl $n$-manifold ($n>2$) with Weyl derivative
$D$, $F^D\wedge\Omega_J$ and the commutator $[R^{D,w}_{X,Y},J]$ both vanish.
If $n>4$ it follows that $F^D=0$, while for $n=4$,
$F^D$ is orthogonal to $\Omega_J$. Also if $R^D=R^{D,1}$ then
the symmetric Ricci tensor is given by the formula
$$\tfrac12\cip{R^D_{Je_i,e_i}X,JY}=(n-2)r^D_0(X,Y)+\tfrac1n\scal^D\cip{X,Y},$$
where on the left we are summing over a weightless orthonormal basis $e_i$.
Consequently the symmetric Ricci tensor is $J$-invariant.
\proofof{prop} The first two facts are immediate from $d^D\Omega_J=0$ and
$DJ=0$ respectively. If $n>4$ then wedge product with $\Omega_J$ is injective
on $2$-forms, while for $n=4$, $F^D\wedge\Omega_J$ is the multiple
$\pm\cip{F^D,\Omega_J}$ of the weightless volume form, since $\Omega_J$ is
antiselfdual. The final formula is a consequence of the first Bianchi
identity:
\begin{align*}
\tfrac12\cip{R^D_{Je_i,e_i}X,JY}
&=\cip{R^D_{X,e_i}Je_i,JY}=\cip{R^D_{X,e_i}e_i,Y}\\
&=F^D(X,e_i)\cip{e_i,Y}-\cip{r^D(X)\skwend e_i\,e_i,Y}
+\cip{r^D(e_i)\skwend X e_i,Y}\\
&=(n-2)r^D_0(X,Y)+\tfrac1n\scal^D\cip{X,Y}-\tfrac12(n-4)F^D(X,Y)
\end{align*}
and the last term vanishes since $F^D=0$ for $n>4$.
\end{proof}
Now suppose $n=4$. Then $W^\sd_{X,Y}$ commutes with $J$, and so
$$J\circ W^\asd_{X,Y}-W^\asd_{X,Y}\circ J =J\circ\bigl(r^D(X)\skwend Y-
r^D(Y)\skwend X\bigr)-\bigl(r^D(X)\skwend Y-r^D(Y)\skwend X\bigr)\circ J.$$
The bundle of antiselfdual Weyl tensors may be identified with the rank $5$
bundle of symmetric trace-free maps $L^2\Lambda^2_\asd
T\dual\Mn\to\Lambda^2_\asd T\dual\Mn$, where $W^\asd(U\wedge V)(X\wedge
Y)=\cip{W^\asd_{U,V}X,Y}$ and we identify $L^2\Lambda^2_\asd T\dual\Mn$ with
$L^{-2}\Lambda^2_\asd T\Mn$. Under the unitary group $L^2\Lambda^2_\asd
T\dual\Mn$ decomposes into the span of $J$ and the weightless canonical bundle
$L^2K_J$. This bundle of Weyl tensors therefore decomposes into three
pieces: the Weyl tensors acting by scalars on $\langle J\rangle$ and $L^2K_J$;
the symmetric trace-free maps $L^2K_J\to K_J$ (acting trivially on $\langle
J\rangle$); and the Weyl tensors mapping $\langle J\rangle$ into $K_J$ and
vice versa. These subbundles have ranks $1,2$ and $2$ respectively. Since no
nonzero Weyl tensor acts trivially on $K_J$, it follows that the above formula
determines $W^\asd$ uniquely in terms of $r^D$. Now this is an invariant
formula which is linear in $r^D$, so $r^D_0$ and $F^D_+$ cannot contribute:
they are sections of (isomorphic) irreducible rank 3 bundles. Thus the first
and third components of $W^\asd$ are given by $\scal^D$ and $F^D_\asd$
respectively, and the second component must vanish. The numerical factors can
now be found by taking a trace.
\begin{prop}\label{Wform}\tcite{AG}
On a K\"ahler-Weyl $4$-manifold with Weyl derivative $D$,
$$W^\asd=\tfrac14\scal^D
\bigl(\tfrac13\iden_{\Lambda^2_\asd}-\tfrac12\Omega_J\tens\Omega_J\bigr)
-\tfrac12(JF^D_\asd\tens\Omega_J+\Omega_J\tens JF^D_\asd),$$ where
$JF^D_\asd=F^D_\asd\circ J$. In particular $W^\asd=0$ iff $F^D_\asd=0$
and $\scal^D=0$.
\end{prop}

The \emphdef{Ricci form} $\rho^D$ on $\Mn$ is defined to be the curvature
of $D$ on the weightless canonical bundle $L^2K_J$. Therefore
\begin{align*}
\rho^D(X,Y)&=-\tfrac i2\cip{R^D_{X,Y}e_k,Je_k}\\
&=-\tfrac i2\bigl(\cip{R^D_{X,e_k}e_k,JY}-\cip{R^D_{Y,e_k}e_k,JX}\bigr)\\
&=i\bigl(2r^D_0(JX,Y)+\tfrac14\scal^D\cip{JX,Y}+2F^D_\asd(JX,Y)\bigr).
\end{align*}
Thus $W^\asd=0$ iff $\rho^D$ and $F^D$ are selfdual $2$-forms.

\section{Shear-free geodesic congruences and Einstein-Weyl
geometry}\label{SGCEW}

On a conformal manifold, a foliation with oriented one dimensional leaves may
be described by a weightless unit vector field $\chi$.  (If $K$ is any
nonvanishing vector field tangent to the leaves, then $\chi=\pm K/|K|$.) Such
a foliation, or equivalently, such a $\chi$, is often called a
\emphdef{congruence}.

If $D$ is any Weyl derivative, then $D\chi$ is a section of $T\dual\Mn\tens
L^{-1}T\Mn$ satisfying $\cip{D\chi,\chi}=0$, since $\chi$ has unit length. Let
$\chi^\perp$ be the orthogonal complement of $\chi$ in $L^{-1}T\Mn$. Under the
orthogonal group of $\chi^\perp$ acting trivially on the span of $\chi$, the
bundle $T\dual\Mn\tens\chi^\perp$ decomposes into four irreducible
components:\quad $L^{-1}\Lambda^2(\chi^\perp)$,\quad
$L^{-1}S^2_0(\chi^\perp)$,\quad $L^{-1}$ (multiples of the identity
$\chi^\perp\to\chi^\perp$),\quad and \ $L^{-1}\chi^\perp$ (the
$\chi^\perp$-valued $1$-forms vanishing on vectors orthogonal to $\chi$).

The first three components of $D\chi$ may be found by taking the skew,
symmetric trace-free and tracelike parts of $D\chi-\chi\tens D_\chi\chi$,
while the final component is simply $D_\chi\chi$. These components are
respectively called the \emphdef{twist, shear, divergence, and acceleration}
of $\chi$ with respect to $D$. If any of these vanish, then the
congruence $\chi$ is said to be \emphdef{twist-free, shear-free,
divergence-free, or geodesic} accordingly.

\begin{prop} Let $\chi$ be a unit section of $L^{-1}T\Mn$. Then the
shear and twist of $\chi$ are independent of the choice of Weyl derivative
$D$. Furthermore there is a unique Weyl derivative $D\S/\chi$ with respect
to which $\chi$ is divergence-free and geodesic.
\end{prop}
\noindent This follows from the fact that
$(D+\gamma)\chi=D\chi+\gamma(\chi)\iden-\chi\tens\gamma$.

The twist is simply the Frobenius tensor of $\chi^\perp$ (i.e., the $\chi$
component of the Lie bracket of sections of $\chi^\perp$), while the shear
measures the Lie derivative of the conformal structure of $\chi^\perp$ along
$\chi$ (which makes sense even though $\chi$ is weightless).

\begin{rem}\label{Exact} If $D\S/\chi$ is exact, with $D\S/\chi\mu=0$
then $K=\mu\chi$ is a geodesic divergence-free vector field of unit length
with respect to the metric $g=\mu^{-2}\conf$. If $\chi$ is also shear-free,
then $K$ is a Killing field of $g$. Note conversely that any nonvanishing
conformal vector field $K$ is a Killing field of constant length $a$ for the
compatible metric $a^2|K|^{-2}\conf$: $\chi=K/|K|$ is then a shear-free
congruence, and $D\S/\chi$ is the exact Weyl derivative $D^{|K|}$, which we
call the \emphdef{constant length gauge} of $K$.
\end{rem}

We now turn to the study of geodesic congruences in three dimensional Weyl
spaces and their relationship to Einstein-Weyl geometry and minitwistor theory
(see~\cite{Hitchin3,LeBrun2,PT1}). We discuss the ``mini-Kerr theorem'' which
is rather a folk theorem in the existing literature, and rewrite the
Einstein-Weyl condition in a novel way by finding special monopole equations
associated to a shear-free geodesic congruence.

The \emphdef{minitwistor space} of an oriented geodesically convex Weyl space
is its space of oriented geodesics. We assume that this is a manifold (i.e.,
we ignore the fact that it may not be Hausdorff), as we shall only be using
minitwistor theory to probe the local geometry of the Weyl space. The
minitwistor space is four dimensional, and has a distinguished family of
embedded $2$-spheres corresponding to the geodesics passing through given
points in the Weyl space.

Now let $\chi$ be a geodesic congruence on an oriented Weyl space $B$ with
Weyl connection $D^B$. Then
\begin{equation}\label{Cong}
D^B\chi=\dv(\iden-\chi\tens\chi)+\tw\,{*\chi}+\Sh,
\end{equation}
where the divergence and twist, $\dv$ and $\tw$, are sections of $L^{-1}$
and $\Sh$ is the shear. Note that $D\S/\chi=D^B-\dv\chi$.

Equation~\eqref{Cong} admits a natural complex interpretation, which we give
in order to compare our formulae to those in the literature~\cite{HT,PT1}. Let
$\cH=\chi^\perp\tens\C$ in the complexified weightless tangent bundle. Then
$\cH$ has a complex bilinear inner product on each fibre and the orientation
of $B$ distinguishes one of the two null lines: if $e_1,e_2$ is an oriented
real orthonormal basis, then $e_1+ie_2$ is null. Let $Z$ be a section of this
null line with $\cip{Z,\overline Z}=1$. Such a $Z$ is unique up to pointwise
multiplication by a unit complex number: at each point it is of the form
$(e_1+ie_2)/\sqrt2$. Now $D^B\chi=\dr\, Z\tens \overline Z
+\overline\dr\,\overline Z\tens Z +\sh\,Z\tens Z +\overline\sh\,\overline
Z\tens\overline Z$, where $\dr=\dv+i\tw$ and $\sh=\Sh(\overline Z,\overline
Z)$ are sections of $L^{-1}\tens\C$. Note that $\sh$ depends on the choice of
$Z$: the ambiguity can partially be removed by requiring that $D^B_\chi Z=0$,
but we shall instead work directly with $\Sh$.

\begin{conv} There are two interesting sign conventions for the Hodge
star operator of an oriented conformal manifold. The first satisfies
$\alpha\wedge{\tilde*\beta}=\cip{\alpha,\beta}\ortn$, where $\ortn$ is the
unit section of $L^n\Lambda^nT\dual\Mn$ given by the orientation.  This is
convenient when computing the star operator of an explicit example. The
second satisfies $*1=\ortn$ and $\iota\s/X\,{*\alpha}=*(X\wedge\alpha)$,
which is a more useful property in many theoretical calculations. Also
$*^2=(-1)^{\frac12n(n-1)}$ depends only on the dimension of the manifold,
not on the degree of the form. If $\alpha$ is a $k$-form, then
$*\alpha=(-1)^{\frac12k(k-1)}\tilde*\alpha$.
\end{conv}
\begin{prop}\label{SFC} The curvature of $D^B$ applied to the geodesic 
congruence $\chi$ is given by
\begin{align*}
R^{B,0}_{X,Y}\chi=\iota\s/\chi\bigl[&
D^B_X\dv\;\chi\wedge Y-D^B_Y\dv\;\chi\wedge X
-D^B_X\tw\;{*Y}+D^B_Y\tw\;{*X}\\
&\qquad\qquad\qquad\quad+(\dv^2-\tw^2)X\wedge Y
-2\dv\tw\,\chi\wedge{*(X\wedge Y)}\bigr]\\
+&(D^B_X\Sh)(Y)-(D^B_Y\Sh)(X)\\
-&\dv\bigl(\Sh(X)\cip{\chi,Y}-\Sh(Y)\cip{\chi,X}\bigr)+
\tw\,{*\bigl(Y\wedge\Sh(X)-X\wedge\Sh(Y)\bigr)}
\end{align*}
and also by its decomposition:
\begin{align*}
R^{B,0}_{X,Y}\chi=\quad&r^B_0(Y,\chi)X-\tfrac12F^B(Y,\chi)X
+\bigl(r^B_0(X)+\tfrac16\scal^BX-\tfrac12F^B(X)\bigr)\cip{\chi,Y}\\
-&r^B_0(X,\chi)Y+\tfrac12F^B(X,\chi)Y
-\bigl(r^B_0(Y)+\tfrac16\scal^BY-\tfrac12F^B(Y)\bigr)\cip{\chi,X}.
\end{align*}
\end{prop}
The first formula is obtained from
$R^{B,0}_{X,Y}\chi=D^B_X(D^B\chi)\s/Y-D^B_Y(D^B\chi)\s/X$,
using $D^B_X(D^B\chi)=D^B_X\dv(\iden-\chi\tens\chi)+D^B_X\tw\;{*\chi}
-\dv(D^B_X\chi\tens\chi+\chi\tens D^B_X\chi)+\tw\,{*D^B_X\chi}+D^B_X\Sh$.
The second formula follows easily from
$R^{B,0}_{X,Y}=-r^B(X)\skwend Y+r^B(Y)\skwend X$ where
$r^B=r^B_0+\frac1{12}\scal^B-\frac12F^B$.

In order to compare the rather different formulae in
Proposition~\ref{SFC}, we shall first take $Y$ parallel to $\chi$ and $X$
orthogonal to $\chi$. The formulae reduce to
\begin{align*}
D^B_\chi\dv&\;X+D^B_\chi\tw\;JX+(\dv^2-\tw^2)X+2\dv\tw\,JX\\
&\qquad\qquad\qquad\qquad
+\Sh(D^B_X\chi)+(D^B_\chi\Sh)(X)+\dv\Sh(X)-\tw\Sh(JX)\\
&=-R^{B,0}_{X,\chi}\chi\\
&=-r^B_0(X)+r^B_0(X,\chi)\chi+\tfrac12\bigl(F^B(X)-F^B(X,\chi)\chi\bigr)
-r^B_0(\chi,\chi)X-\tfrac16\scal^BX,
\end{align*}
where $JX:=\iota\s/X\,{*\chi}$ and we have used the fact that
$(D^B_X\Sh)(\chi)+\Sh(D^B_X\chi)=0$. If we contract with another
vector field $Y$ orthogonal to $\chi$, then we obtain
\begin{multline*}
D^B_\chi\dv\,\cip{X,Y}+D^B_\chi\tw\,\cip{JX,Y}
+\cip{(D^B_\chi\Sh)(X),Y}\\
+(\dv^2-\tw^2)\cip{X,Y}
+2\dv\tw\cip{JX,Y}
+2\dv\cip{\Sh(X),Y}+\cip{\Sh(X),\Sh(Y)}\\
=-r^B_0(X,Y)+\tfrac12F^B(X,Y)-\bigl(r^B_0(\chi,\chi)-\tfrac16\scal^B\bigr)
\cip{X,Y}.
\end{multline*}
Decomposing this into irreducibles gives the equations
\begin{align}\label{Jacobi1}
D^B_\chi\dv+\dv^2-\tw^2+\tfrac12|\Sh|^2
+\tfrac12r^B_0(\chi,\chi)+\tfrac16\scal^B&=0\\
D^B_\chi\tw+2\dv\tw+\tfrac12\cip{\chi,*F^B}&=0\\
D^B_\chi\Sh+2\dv\Sh+\sym_0^{\chi^\perp}r^B_0&=0\label{Jacobi3}\\
\intertext{which may, assuming $D^B_\chi Z=0$, be rewritten as}
\label{Jacobi4}
D^B_\chi\dr+\dr^2+\sh\overline\sh
+\tfrac12r^B_0(\chi,\chi)+\tfrac16\scal^B+\tfrac i2\cip{\chi,*F^B}&=0\\
D^B_\chi\sh+(\dr+\overline\dr)\sh+r^B_0(\overline Z,\overline Z)&=0.
\label{Jacobi5}
\end{align}

Along a single geodesic, these formulae describe the evolution of nearby
geodesics in the congruence and therefore may be interpreted infinitesimally
(cf.~\cite{PT1}). We say that a vector field $X$ along an oriented geodesic
$\Gamma$ with weightless unit tangent $\chi$ is a \emphdef{Jacobi field} iff
$(D^B)^2_{\chi,\chi}X=R^{B,0}_{\chi,X}\chi$. The space of Jacobi fields
orthogonal to $\Gamma$ is four dimensional, since the initial data for the
Jacobi field equation is $X,D^B_\chi X$. In fact this is the tangent space to
the minitwistor space at $\Gamma$.  If we now consider a two dimensional
family of Jacobi fields spanning (at each point on an open subset of $\Gamma$)
the plane orthogonal to $\Gamma$, then we may write $D^B_\chi X=\dv
X+\tw\,JX+\Sh(X)$ for the Jacobi fields $X$ in this family. If we
differentiate again with respect to $\chi$, we reobtain the
equations~(\ref{Jacobi1})--(\ref{Jacobi5}).

A geodesic congruence gives rise to such a two dimensional family of Jacobi
fields along each geodesic in the congruence. We define the Lie derivative
$\cL_\chi X$ of a vector field $X$ along $\chi$ to be the horizontal part of
$D^B_\chi X-D^B_X\chi$. Then if $\cL_\chi X=0$, $X$ is a Jacobi field, and
such Jacobi fields are determined along a geodesic by their value at a
point. Next note that $\cL_\chi J=0$ (i.e., $\cL_\chi(JX)=J\cL_\chi X$) iff
$\chi$ is shear-free. However, equation~(\ref{Jacobi3}) shows that if $\chi$
is a shear-free, then $r^B_0(X,Y)=-\frac12r^B_0(\chi,\chi)\cip{X,Y}$ for all
$X,Y$ orthogonal to $\chi$. More generally, this equation shows that $J$ is
a well defined complex structure on the space of Jacobi fields orthogonal to
a geodesic $\Gamma$ iff $r^B_0(X,Y)=-\frac12r^B_0(\chi,\chi)\cip{X,Y}$ for
all $X,Y$ orthogonal to $\Gamma$. The Jacobi fields defined by a congruence
are then invariant under $J$ iff the congruence is shear-free.

\begin{defn}\label{EWdef}\tcite{Hitchin3} A Weyl space $B,D^B$ is said to be
\emphdef{Einstein-Weyl} iff $r^B_0=0$.
\end{defn}

As mentioned above, the space of orthogonal Jacobi fields along a geodesic is
the tangent space to the minitwistor space at that geodesic. Therefore, if $B$
is Einstein-Weyl, the minitwistor space admits a natural almost complex
structure. This complex structure turns out to be integrable, and so the
minitwistor space of an Einstein-Weyl space is a complex surface $\cS$
containing a family of rational curves, called \emphdef{minitwistor lines},
parameterised by points in $B$~\cite{Hitchin3}. These curves have normal
bundle $\cO(2)$ and are invariant under the real structure on $\cS$ defined by
reversing the orientation of a geodesic. Conversely, any complex surface with
real structure, containing a real (i.e., invariant) rational curve with normal
bundle $\cO(2)$, determines an Einstein-Weyl space as the real points in the
Kodaira moduli space of deformations of this curve. We therefore have a
twistor construction for Einstein-Weyl spaces, called the \emphdef{Hitchin
correspondence}. We note that the canonical bundle $K_\cS$ of $\cS$ has degree
$-4$ on each minitwistor line.

Since geodesics correspond to points in the minitwistor space, a geodesic
congruence defines a real surface $\cC$ intersecting each minitwistor line
once. By the definition of the complex structure on $\cS$, the surface $\cC$ is
a holomorphic curve iff the geodesic congruence is shear-free. This may be
viewed as a minitwistor version of the Kerr theorem: \emph{every
shear-free geodesic congruence in an Einstein-Weyl space is obtained locally
from a holomorphic curve in the minitwistor space}. In particular, we have
the following.
\begin{prop} Let $B,D^B$ be a three dimensional Weyl space. Then the
following are equivalent:
\begin{enumerate}
\item $B$ is Einstein-Weyl 
\item Given any point $b\in B$ and any unit vector $v\in L^{-1}T_bB$,
there is a shear-free geodesic congruence $\chi$ defined on
a neighbourhood of $b$ with $\chi_b=v$
\item Given any point $b\in B$ there are three shear-free geodesic
congruences defined on a neighbourhood of $b$ which are pairwise
non-tangential at $b$.
\end{enumerate}
\proofof{prop} Clearly (ii) implies (iii). It is immediate
from~\eqref{Jacobi3} that (ii) implies (i); to obtain the stronger result that
(iii) implies (i) suppose that $B$ is not Einstein-Weyl, i.e., at some point
$b\in B$, $r^B_0\neq0$. If $\chi$ is a shear-free geodesic congruence near $b$
then by equation~\eqref{Jacobi3}, $r^B_0$ is a multiple of the identity on
$\chi^\perp$, and one easily sees that this multiple must be the middle
eigenvalue $\lambda_0\in L^{-2}_b$ of $r^B_0$ at $b$. Now at $b$, $r^B_0$ may
be written $\alpha\tens\sharp\beta+\beta\tens\sharp\alpha+\lambda_0\iden$
where $\alpha,\beta\in T\dual_b B$ with
$\cip{\alpha,\beta}=-\frac32\lambda_0$. The directions of $\sharp\alpha$ and
$\sharp\beta$ are uniquely determined by $r^B_0$ and $\chi$ must lie in one of
these directions. Hence if $B$ is not Einstein-Weyl at $b$, there are at most
two possible directions at $b$ (up to sign) for a shear-free geodesic
congruence. (Note that the linear algebra involved here is the same as that
used to show that there are at most two principal directions of a nonzero
antiselfdual Weyl tensor in four dimensions; see, for instance~\cite{AG}.
Our result is just the symmetry reduction of this fact.)

Finally, to see that (i) implies (ii), we simply observe that given any
minitwistor line and any point on that line, we can find, in a neighbourhood
of that point, a transverse holomorphic curve. This curve will also intersect
nearby minitwistor lines exactly once.
\end{proof}

We now want to study Einstein-Weyl spaces with a shear-free geodesic
congruence in more detail. As motivation for our main result, notice that the
curve $\cC$ in the minitwistor space given by $\chi$ defines divisors
$\cC+\overline\cC$ and $\cC-\overline\cC$ such that the line bundles
$[\cC+\overline\cC]K^{1/2}_\cS$ and $[\cC-\overline\cC]$ are trivial on each
minitwistor line. It is well known~\cite{JT} that such line bundles correspond
to solutions $(\gmw,A)$ of the abelian monopole equation $*D^B\gmw=dA$, where
$\gmw$ is a section of $L^{-1}$ and $A$ is a $1$-form. Therefore, we should be
able to find two special solutions of this monopole equation, one real and one
imaginary, associated to any shear-free geodesic congruence.

These solutions turn out to be $\tw$ and $i\dv$. To see this, we return to
the curvature equations in Proposition~\ref{SFC} and look at the horizontal
components. If $X,Y$ are orthogonal to a geodesic congruence $\chi$ on any
three dimensional Weyl space then:
\begin{multline*}
D^B_X\dv\;Y-D^B_Y\dv\;X+D^B_X\tw\;JY-D^B_Y\tw\;JX\\
+(D^B_X\Sh)(Y)-(D^B_Y\Sh)(X)+\tw\,{*(Y\wedge\Sh(X)-X\wedge\Sh(Y))}\\
=r^B_0(Y,\chi)X-\tfrac12F^B(Y,\chi)X-r^B_0(X,\chi)Y+\tfrac12F^B(X,\chi)Y.
\end{multline*}
If $\chi$ is shear-free this reduces to the equation
$$D^B_X\dv-D^B_{JX}\tw+r^B_0(\chi,X)+\tfrac12F^B(\chi,X)=0,$$
where $X\perp\chi$. From this, and our earlier formulae, we have:

\begin{prop} Let $\chi$ be shear-free geodesic congruence with
divergence $\dv$ and twist $\tw$ in a three dimensional Weyl
space $B$. Then $\chi$ satisfies the equations
\begin{align}
D^B_\chi\dv+\dv^2-\tw^2+\tfrac16\scal^B&=0\label{short1}\\
D^B_\chi\tw+2\dv\tw+\tfrac12\cip{\chi,*F^B}&=0\label{short2}\\
(D^B\dv-D^B\tw\circ J)|\s/{\chi^\perp}+\tfrac12\iota\s/\chi F^B&=0
\label{CReqn}
\end{align}
if and only if $B$ is Einstein-Weyl.
\end{prop}
The last equation, like the first two (see~\eqref{Jacobi4}),
admits a natural complex formulation in terms of $\dr$. Instead, however,
we shall combine these equations to give the following result.

\begin{thm}\label{EWeqnthm}
The three dimensional Einstein-Weyl equations are equivalent
to the following special monopole equations for a shear-free geodesic
congruence $\chi$ with $D^B\chi=\dv(\iden-\chi\tens\chi)+\tw\,{*\chi}$.
\begin{align}
*D^B\dv&=-\tfrac12{*\iota\s/\chi F^B}-\tfrac16\scal^B{*\chi}
-(\dv^2+\tw^2){*\chi}+d(\tw\chi)\label{long1}\\
*D^B\tw&=\tfrac12F^B-d(\dv\chi).\label{long2}
\end{align}
\textup[By ``monopole equations'', we mean that the right hand sides are
closed $2$-forms. Note also that these equations are not independent:
they are immediately equivalent to~\eqref{short1} and~\eqref{long2}, or
to~\eqref{short2} and~\eqref{long1}.\textup]
\proofof{thm} The equations of the previous proposition are equivalent to the
following:
\begin{align*}
D^B\dv&=D^B\tw\circ J+(\tw^2-\dv^2)\chi
-\tfrac16\scal^B\chi-\tfrac12\iota\s/\chi F^B
\\
D^B\tw&=-D^B\dv\circ J-2\dv\tw\chi-\tfrac12{*F^B}.
\end{align*}
Applying the star operator readily yields the equations of the theorem. The
second equation is clearly a monopole equation, since $F^B$ is closed. It
remains to check that the right hand side of the first equation is closed:
\begin{align*}
d\bigl(&\tfrac12\chi\wedge{*F^B}+\tfrac16\scal^B{*\chi}
+(\dv^2+\tw^2){*\chi}\bigr)\\
&=\tfrac12d^B\chi\wedge{*F^B}-\tfrac12\chi\wedge{*\delta^BF^B}
+\tfrac16D^B\scal^B\wedge{*\chi}+(2\dv D^B\dv+2\tw D^B\tw)\wedge{*\chi}\\
&\qquad+(\tfrac16\scal^B+\dv^2+\tw^2){*\delta^B\chi}\\
&=\tfrac12\chi\wedge{*\bigl(\tfrac13D^B\scal^B-\delta^BF^B\bigr)}\\
&\qquad+\bigl(\tw\cip{\chi,*F^B}+2\dv D^B_\chi\dv+2\tw D^B_\chi\tw
+2\dv(\tfrac16\scal^B+\dv^2+\tw^2)\bigr){*1}.
\end{align*}
Here $\delta^B=\trace D^B$ is the divergence on forms, and so the first term
vanishes by virtue of the second Bianchi identity. The remaining multiple of
the orientation form $*1$ is
$$\tw\cip{\chi,*F^B}+2\tw D^B_\chi\tw+2\tw(2\tw\dv)
+2\dv D^B_\chi\dv+2\dv(\tfrac16\scal^B+\dv^2-\tw^2),$$
which vanishes by the previous proposition.
\end{proof}

Two key special cases of this theorem have already been studied.

\subsection*{LeBrun-Ward geometries}\noindent\medbreak

Suppose an Einstein-Weyl space admits a shear-free geodesic congruence
\emph{which is also twist-free}. Then $\tw=0$ and so the Einstein-Weyl
equations~\eqref{short1},~\eqref{long2} are:
\begin{align}
D^B_\chi\dv+\dv^2&=-\tfrac16\scal^B\label{abstractTFE}\\
F^B=2d(\dv\chi)&=2D^B\dv\wedge\chi.\label{LWgauge}
\end{align}
As observed by Tod~\cite{Tod2}, these Einstein-Weyl spaces are the spaces
first studied by LeBrun~\cite{LeBrun1,LeBrun2} and Ward~\cite{Ward}, who
described them using coordinates in which the above equations reduce to the
$\SU(\infty)$ Toda field equation $u_{xx}+u_{yy}+(e^u)_{zz}=0$. Consequently
these Einstein-Weyl spaces are also said to be \emphdef{Toda}.

It may be useful here to sketch how this follows from our formulae, since
Lemma 4.1 in~\cite{Tod2}, given there without proof, is only true after making
use of the gauge freedom to set $z=f(\tilde z)$ and rescale the metric by
$f'(\tilde z)^{-2}$. The key point is that since $D^{LW}:=D^B-2\dv\chi$ is
locally exact by~\eqref{LWgauge}, there is locally a canonical gauge (up to
homothety) in which to work, which we call the \emphdef{LeBrun-Ward} gauge
$\mu_{LW}$. Since $\chi$ is twist-free and also geodesic with respect to
$D^{LW}$, the $1$-form $\mu_{LW}^{-1}\chi$ is locally exact. Taking this to be
$dz$ and introducing isothermal coordinates $(x,y)$ on the quotient of $B$ by
$\chi$, we may write $g\s/{LW}=e^u(dx^2+dy^2)+dz^2$ for some function
$u(x,y,z)$, since $\chi$ is shear-free. By computing the divergence of $\chi$
we then find that the Toda monopole is $\dv=-\frac12u_z\mu_{LW}^{-1}$, and
equation~\eqref{abstractTFE} reduces easily in this gauge to the Toda
equation. One of the reasons for the interest in this equation is that it may
be used to construct hyperK\"ahler and scalar-flat K\"ahler
$4$-manifolds~\cite{BF,LeBrun1}, as we shall see in the next section.

LeBrun~\cite{LeBrun1} shows that these spaces are characterised by the
existence of a divisor $\cC$ in the minitwistor space with
$[\cC+\overline\cC]=K_\cS^{-1/2}$. This agrees with our assertion that
$[\cC+\overline\cC]K_\cS^{1/2}$ corresponds to the monopole $\tw$.

In~\cite{DMJC4}, it is shown that an Einstein-Weyl space admits at
most a three dimensional family of shear-free twist-free geodesic
congruences.

\subsection*{Gauduchon-Tod geometries}\noindent\medbreak

Suppose an Einstein-Weyl space admits a shear-free geodesic
congruence \emph{which is also divergence-free}. Then $\dv=0$ and so
the Einstein-Weyl equations~\eqref{short1},~\eqref{long2} are:
\begin{align}
\tw^2&=\tfrac16\scal^B\label{scalmon}\\
*D^B\tw&=\tfrac12F^B.\label{GTeqn}
\end{align}
It follows that these are the geometries which arose in the work of Gauduchon
and Tod~\cite{GT} and also Chave, Tod and Valent~\cite{CTV} on hypercomplex
$4$-manifolds with triholomorphic conformal vector fields. Gauduchon and Tod
essentially observe the following equivalent formulation of these equations.

\begin{prop} The connection $D^\tw=D^B-\tw\,{*1}$ on
$L^{-1}TB$ is flat.
\proofof{prop} The curvature of $D^\tw$ is easily
computed to be:
\begin{multline*}
R^\tw_{X,Y}=-r^B_0(X)\skwend Y+r^B_0(Y)\skwend X-\tfrac16\scal^BX\skwend Y\\
+\tfrac12F^B(X)\skwend Y-\tfrac12F^B(Y)\skwend X
-D^B_X\tw\,{*Y}+D^B_Y\tw\,{*X}+\tw^2 X\skwend Y.
\end{multline*} Now
$D^B_X\tw\,{*Y}-D^B_Y\tw\,{*X}=(*D^B\tw)(X)\skwend Y-(*D^B\tw)(Y)\skwend Y$,
so equations~\eqref{scalmon} and~\eqref{GTeqn} imply that $R^\tw_{X,Y}$
vanishes if $B$ is Einstein-Weyl.  [Conversely if there is a $\chi$ with
$R^\tw_{X,Y}\chi=0$ for all $X,Y$, then $B$ is Einstein-Weyl.]
\end{proof}

This shows that the existence of a single shear-free divergence-free geodesic
congruence gives an entire $2$-sphere of such congruences and we say that
these Einstein-Weyl spaces are \emphdef{hyperCR}~\cite{CT}. There is also a
simple minitwistor interpretation of this. The divisor $\cC$ corresponding to
a shear-free divergence-free geodesic congruence has $[\cC-\overline\cC]$
trivial, i.e., $\cC-\overline\cC$ is the divisor of a meromorphic function.
Hence we have a nonconstant holomorphic map from the minitwistor space to
$\CP1$, and its fibres correspond to the $2$-sphere of congruences. This
argument is the minitwistor analogue of the twistor characterisation of
hypercomplex structures discussed in the previous section.

Since the Einstein-Weyl structure determines $\tw$ up to sign, it follows that
an Einstein-Weyl space admits at most two hyperCR structures. If it admits
exactly two, then we must have $\tw\neq0$ and $F^B=0$, i.e., the Einstein-Weyl
space is the round sphere.

\subsection*{Einstein-Weyl spaces with a geodesic symmetry}\noindent\medbreak

The Einstein-Weyl equation can be completely solved in the case of
Einstein-Weyl spaces admitting a shear-free geodesic congruence $\chi$ such
that $\chi=K/|K|$ with $K$ a conformal vector field preserving the Weyl
connection. In this case $D\S/\chi=D^B-\dv\chi$ is exact, $|K|$ being a
parallel section of $L^1$ (see Remark~\ref{Exact}). We introduce
$g=|K|^{-2}\conf$ so that $D\S/\chi=D^g$. Since $K$ preserves the Weyl
connection and $\cL_Kg=0$, we may write $\dv=\dv_g|K|^{-1},\tw=\tw_g|K|^{-1}$,
where $\partial_K\dv_g=\partial_K\tw_g=0$. Now $\iota\s/\chi F^B=\iota\s/\chi
d(\dv\chi) =-D^g\dv$ and so equation~\eqref{CReqn} becomes
$$\tfrac12d\dv_g-d\tw_g\circ J=0.$$
This is solved by setting $2\tw_g-i\dv_g=H$, where $H$ is a holomorphic
function on the quotient $C$ of $B$ by $K$. Since $D^B_\chi\dv=-\dv^2$ and
$D^B_\chi\tw=-\dv\tw$, the remaining Einstein-Weyl equations reduce to
$\dv\tw+\frac12\cip{\chi,*F^B}=0$ and $\tw^2=\frac16\scal^B$. The first of
these is automatic. To solve the second we note that $\scal^B$ can be computed
from the scalar curvature of the quotient metric on $C$ using a submersion
formula~\cite{Besse,CP}. This gives $\scal^B=\scal^C-2\dv^2-2\tw^2$ and hence
$\scal^C=2\dv^2+8\tw^2=2|2\tw-i\dv|^2$. If this is zero, then $\dv=\tw=0$ and
$D^B$ is flat. Otherwise we observe that $\log|H|^2$ is harmonic, and so
rescaling the quotient metric by $|H|^2$ gives a metric of constant curvature
$1$ (i.e., the scalar curvature is $2$).

Remarkably, these Einstein-Weyl spaces are also all hyperCR: since
$\tw^2=\frac16\scal^B$ and $*D^B\tw=\frac12F^B-d(\dv\chi)=-\frac12F^B$,
reversing the sign of $\tw$ (or equivalently, reversing the orientation of
$B$) solves the equations of the previous subsection. Thus we have established
the following theorem.
\begin{thm}\label{gs}
The three dimensional Einstein-Weyl spaces with geodesic
symmetry are either flat with translational symmetry or are given locally by:
\begin{align*}
g&=|H|^{-2}(\sigma_1^2+\sigma_2^2)+\beta^2\\
\omega&=\tfrac i2(H-\overline H)\beta\\
d\beta&=\tfrac12(H+\overline H)|H|^{-2}\sigma_1\wedge\sigma_2
\end{align*}
where $\sigma_1^2+\sigma_2^2$ is the round metric on $S^2$, and $H$ is any
nonvanishing holomorphic function on an open subset of $S^2$. The geodesic
symmetry $K$ is dual to $\beta$ and the monopoles associated to $K/|K|$ are
$\dv=\frac i2(H-\overline H)\mu_g^{-1}$ and
$\tw=\frac14(H+\overline H)\mu_g^{-1}$. These spaces all admit hyperCR
structures, with flat connection $D^B+\tw\,{*1}$.
\end{thm}
The equation for $\beta$ can be integrated explicitly. Indeed if
$\zeta$ is a holomorphic coordinate such that
$\sigma_1\wedge\sigma_2=2i\,d\zeta\wedge d\overline\zeta/
(1+\zeta\overline\zeta)^2$ then one can take
$$\beta=d\psi+\frac i{1+\zeta\overline\zeta}\left(\frac{d\zeta}{\zeta H}
-\frac{d\overline\zeta}{\overline\zeta\mskip2mu\overline H}\right).$$
Of course, this is not the only possible choice: for instance one
can write $d\zeta/(\zeta H)=dF$ with $F$ holomorphic and use
$\beta=d\psi-i(F-\overline F)\,d\bigl(1/(1+\zeta\overline\zeta)\bigr)$.

Note that $\omega$ is dual to a Killing field of $g$ iff $H$ is constant, in
which case we obtain the well known Einstein-Weyl structures on the Berger
spheres. The Einstein metric on $S^3$ arises when $H$ is real, in which case
the connections $D^B\pm\tw\,{*1}$ are both flat: they are the left and right
invariant connections. The flat Weyl structure with radial symmetry (which
is globally defined on $S^1\times S^2$) occurs when $H$ is purely imaginary.
Gauduchon and Tod~\cite{GT} prove that these are the only hyperCR
structures on compact Einstein-Weyl manifolds.

The fact that the Einstein-Weyl spaces with geodesic symmetry are hyperCR may
equally be understood via minitwistor theory. Indeed, any symmetry $K$ (a
conformal vector field preserving the Weyl connection) on a $3$-dimensional
Einstein-Weyl space induces a holomorphic vector field $X$ on the minitwistor
space $\cS$. If $K$ is nonvanishing, then on each minitwistor line, $X$ will
be tangent at two points (since the normal bundle is $\cO(2)$) and if the line
corresponds to a real point $x$, then these two tangent points in $\cS$ will
correspond to the two orientations of the geodesic generated by $K_x$. Hence
$X$ vanishes at a point of $\cS$ iff $K$ is tangent along the corresponding
geodesic.

Now if $K$ is a geodesic symmetry then $X$ will be tangent to each minitwistor
line precisely at the points at which it vanishes, and the zeroset of $X$ will
be a divisor (rather than isolated points). This means that $X$ is a section
of a line subbundle $\cH=[\mathrm{div} X]$ of $T\cS$ transverse to the
minitwistor lines ($\cH$ must be transverse even where $X$ vanishes, because
$K$, being real, is not null, and so the points of tangency are simple): the
$\tw$ monopole of $K$ is therefore $\cH\tens K_\cS^{1/2}$. Now the integral
curves of the distribution $\cH$ in the neighbourhood $U$ of some real
minitwistor line give a holomorphic map from $U$ to $\CP1$. Viewing this as a
meromorphic function (by choosing conjugate points on $\CP1$) we obtain a
divisor $\cC-\overline\cC$, where $\cC+\overline\cC$ is a divisor for
$T\cS/\cH$, because $T\cS/\cH$ is isomorphic to $T\CP1$ over each minitwistor
line. Since $K_\cS^{-1}=\cH\tens T\cS/\cH$ we find that
$[\cC+\overline\cC]K_\cS^{1/2}$ is dual to $\cH\tens K_\cS^{1/2}$, which
explains (twistorially) why the $\tw$ monopole of the hyperCR structure is
simply the negation of the $\tw$ monopole of the geodesic symmetry.

Another explanation is that the geodesic symmetry preserves the hyperCR
congruences. Indeed, we have the following observation.

\begin{prop}\label{gschar} Suppose that $B$ is a hyperCR Einstein-Weyl space
with flat connection $D^B+\tw\,{*1}$. Then a vector field $K$ preserves the
hyperCR congruences $\chi$ \textup(i.e., $\cL_K\chi=0$ for each $\chi$\textup)
if and only if it is a geodesic symmetry with twist $\tw$.
\proofof{prop} Since $\chi$ is a weightless vector field,
$\cL_K\chi=D^B_K\chi-D^B_\chi K+\frac13(\trace D^BK)\chi$. This vanishes iff
$D^B_\chi K = \frac13(\trace D^BK)\chi-\tw\,{*(K\wedge\chi)}$.  Hence
$\cL_K\chi=0$ for all of the hyperCR congruences $\chi$ iff
$D^BK=\frac13(\trace D^BK)\iden+\tw\,{*K}$.  This formula shows that $K$ is a
conformal vector field, and that $K/|K|$ is a shear-free geodesic congruence
with twist $\tw$. Also $K$ preserves the flat connection $D^B+\tw\,{*1}$,
since it preserves the parallel sections. Finally, note that the twist of $K$
is determined by the conformal structure from the skew part of $D^{|K|}K$, so
it is also preserved by $K$. Hence $K$ preserves $D^B$ and is therefore a
geodesic symmetry.
\end{proof}

\section{The Jones and Tod construction}\label{JTsect}

In~\cite{JT}, Jones and Tod proved that the quotient of a selfdual conformal
manifold $\Mn$ by a conformal vector field $K$ is Einstein-Weyl: the twistor
lines in the twistor space $Z$ of $\Mn$ project to rational curves with normal
bundle $\cO(2)$ in the space $\cS$ of trajectories of the holomorphic vector
field on $Z$ induced by $K$. Furthermore the Einstein-Weyl space comes with a
solution of the monopole equation from which $\Mn$ can be recovered: indeed
$Z$ is (an open subset of) the total space of the line bundle over $\cS$
determined by this monopole. In other words there is a correspondence
between selfdual spaces with symmetry and Einstein-Weyl spaces with
monopoles. In this section, we explain the differential geometric
constructions involved in the Jones and Tod correspondence, and prove that
invariant antiselfdual complex structures on $\Mn$ correspond to shear-free
geodesic congruences on $B$. These direct methods, although motivated by
the twistor approach, also reveal what happens when $\Mn$ is not selfdual.

Therefore we let $\Mn$ be an oriented conformal manifold with a conformal
vector field $K$, and (by restricting to an open set if necessary) we assume
$K$ is nowhere vanishing. Let $D^{|K|}$ be the constant length gauge of $K$,
so that $\cip{D^{|K|}K,.}$ is a weightless $2$-form. One can compute $D^{|K|}$
in terms of an arbitrary Weyl derivative $D$ by the formula
$$D^{|K|}=D-\frac{\cip{DK,K}}{\cip{K,K}}=
D-\frac14\frac{(\trace DK)K}{\cip{K,K}}
 +\frac12\frac{ (d^DK)(K,.)}{\cip{K,K}},$$
where $(d^DK)(X,Y)=\cip{D_XK,Y}-\cip{D_YK,X}$.

The key observation for the Jones and Tod construction is that there is a
unique Weyl derivative $D^{sd}$ on $\Mn$ such that $\cip{D^{sd}K,.}$ is a
weightless \emph{selfdual} $2$-form: let $\omega=-(*{d^DK})(K,.)/\cip{K,K}$
(which is independent of $D$) and define
$$D^{sd}=D^{|K|}+\tfrac12\omega=D-\frac14\frac{(\trace DK)K}{\cip{K,K}}
+\frac12\frac{(d^DK)(K,.)-(*{d^DK})(K,.)}{\cip{K,K}}.$$
Since $D$ is arbitrary, we may take $D=D^{sd}$ to get
$(D^{sd}K-*{D^{sd}K})(K,.)=0$ from which it is immediate that
$D^{sd}K=*{D^{sd}K}$ since an antiselfdual $2$-form is uniquely determined by
its contraction with a nonzero vector field. The Weyl derivative $D^{sd}$
plays a central role in the proof that $D^B=D^{|K|}+\omega$ is Einstein-Weyl
on $B$. Notice that the the conformal structure and Weyl derivatives
$D^{|K|}$, $D^{sd}$, $D^B$ do indeed descend to $B$ because $K$ is a Killing
field in the constant length gauge and $\omega$ is a basic $1$-form.  Since
the Lie derivative of Weyl derivatives on $L^1$ is given by
$\cL_KD=\frac1nd\trace DK+F^D(K,.)$, it follows that $F^{sd}(K,.)=F^B(K,.)=0$.

We call $D^B$ the \emphdef{Jones-Tod} Weyl structure on $B$.

\begin{thm}\label{JonesnTod}\tcite{JT} Suppose $\Mn$ is an oriented conformal
$4$-manifold and $K$ a conformal vector field such that $B=\Mn/K$ is a
manifold. Let $D^{|K|}$ be the constant length gauge of $K$ and
$\omega=-2(*D^{|K|}K)(K,.)/\cip{K,K}$. Then the Jones-Tod Weyl structure
$D^B=D^{|K|}+\omega$ is Einstein-Weyl on $B$ if and only if $\Mn$ is selfdual.

Note that $*D^B|K|^{-1}=-{*\omega}|K|^{-1}$ is a closed $2$-form. Conversely,
if $(B,D^B)$ is an Einstein-Weyl $3$-manifold and $\gmw\in\Cinf(B,L^{-1})$ is
a nonvanishing solution of the monopole equation $d{*D^B\gmw}=0$ then there is
a selfdual $4$-manifold $\Mn$ with symmetry over $B$ such that $*D^B\gmw$ is
the curvature of the connection defined by the horizontal distribution.

\proofof{thm} The monopole equation on $B$ is equivalent, via the definition
of $\omega$, to the fact that $D^{sd}$ lies midway between $D^B$ and
$D^{|K|}$. So it remains to show that under this condition, the Einstein-Weyl
equation on $B$ is equivalent to the selfduality of $M$. The space of
antiselfdual Weyl tensors is isomorphic to $S^2_0(K^\perp)$ via the map
sending $W^\asd$ to $W^\asd_{.,K}K$, and so it suffices to show that
$r^B_0=0$ iff $W^\asd_{.,K}K=0$.

Since $D^{sd}$ is basic, as a Weyl connection on $T\Mn$,
$0=(\cL_KD^{sd})_X=R^{sd}_{K,X}+D^{sd}_XD^{sd}K$. Therefore:
$$D^{sd}_XD^{sd}K=W\s/{X,K}+r^{sd}(K)\skwend X-r^{sd}(X)\skwend K.$$
If we now take the antiselfdual part of this equation, contract with $K$ and
$Y$, and take $\cip{X,K}=\cip{Y,K}=0$ then we obtain
$$2\cip{W^\asd_{X,K}K,Y}+r^{sd}(X,Y)\cip{K,K}+r^{sd}(K,K)\cip{X,Y}
+*\bigl(K\wedge r^{sd}(K)\wedge X\wedge Y\bigr)=0.$$
Symmetrising in $X,Y$, we see that $W^\asd=0$ iff the horizontal part
of the symmetric Ricci endomorphism of $D^{sd}$ is a multiple of the
identity. The first submersion formula~\cite{Besse} relates the
Ricci curvature of $D^{|K|}$ on $B$ to the horizontal Ricci curvature
of $D^{|K|}$ on $\Mn$. If we combine this with the fact that
$D^{sd}=D^{|K|}+\frac12\omega$ and $D^B=D^{|K|}+\omega$, then
we find that
\begin{multline*}
\sym\Ric^{D^B}_B(X,Y)=\\
\sym\Ric^{D^{sd}}_\Mn(X,Y)+2\cip{D^{|K|}_XK,D^{|K|}_YK}
+\tfrac12\omega(X)\omega(Y)+\mu\cip{X,Y}
\end{multline*}
for some section $\mu$ of $L^{-2}$. Since $D^{|K|}_KK=0$,
$\omega$ vanishes on the plane spanned by $D^{|K|}K$, and so, by comparing
the lengths of $\omega$ and $D^{|K|}K$, we verify that the trace-free part of
$2\cip{D^{|K|}_XK,D^{|K|}_YK}+\frac12\omega(X)\omega(Y)$ vanishes. Hence
$W^\asd=0$ on $\Mn$ iff $D^B$ is Einstein-Weyl on $B$.
\end{proof}
The inverse construction of $\Mn$ from $B$ can be carried out explicitly by
writing $*\s/B D^B\gmw=dA$ on $U\subset B$, so that the real line bundle $\Mn$
is locally isomorphic to $U\times\R$ with connection $dt+A$, where $t$ is the
fibre coordinate. Then the conformal structure
$\conf\s/\Mn=\pi^*\conf\s/B+\gmw^{-2}(dt+A)^2$ is selfdual and
$K=\partial/\partial t$ is a unit Killing field of the representative metric
$g\s/\Mn=\pi^*\gmw^2\conf\s/B+(dt+A)^2$. Note that $\gmw=\pm|K|^{-1}$ and that
the orientations on $M$ and $B$ are related by
$*(\xi\wedge\alpha)=(*\s/B\alpha)\gmw|K|$ where $\alpha$ is any $1$-form on
$B$ and $\xi=K|K|^{-1}$. This ensures that if $D^B=D^{|K|}+\omega$, then the
equation $-(*\s/B\omega)\gmw={*\s/B D^B\gmw}=dA$ is equivalent to
${*(\xi\wedge\omega)}|K|^{-1}=-d(dt+A)$ and hence
$\omega=-\iota\s/K(*d^DK)/|K|^2$ as above.

Jones and Tod also observe that any other solution $(\gmw_1,A_1)$ of the
monopole equation on $B$ corresponds to a selfdual Maxwell field on $\Mn$ with
potential $\tilde A_1=A_1-(\gmw_1/\gmw)(dt+A)$. Indeed, since
$(dt+A)=|K|^{-1}\xi$, one readily verifies that
$$d\tilde A_1
=\bigl(\gmw^{-1}|K|^{-1}\xi\wedge D^B\gmw_1+dA_1\bigr)
-\frac{\gmw_1}{\gmw}\bigl(\gmw^{-1}|K|^{-1}\xi\wedge D^B\gmw+dA\bigr),$$
which is selfdual by the monopole equations for $\gmw$ and $\gmw_1$, together
with the orientation conventions above.

\medbreak

We now want to explain the relationship between invariant complex
structures on $\Mn$ and shear-free geodesic congruences on $B$. That these
should be related is again clear from the twistor point of view: indeed if
$\cD$ is an invariant divisor on $Z$, then it descends to a divisor $\cC$
in $\cS$, which in turn defines, at least locally, a shear-free geodesic
congruence. The line bundles $[\cD+\overline\cD]K_Z^{1/2}$ and
$[\cD-\overline\cD]$ are the pullbacks of $[\cC+\overline\cC]K_\cS^{1/2}$ and
$[\cC-\overline\cC]$ and so we expect the Faraday and Ricci forms on $\Mn$
to be related to the twist and divergence of the congruence on $B$. In
order to see all this in detail, and without the assumption of selfduality,
we carry out the constructions directly.

Suppose that $J$ is an antiselfdual complex structure on $\Mn$ with
$\cL_KJ=0$, so that $K$ is a holomorphic conformal vector field. If $D$ is the
K\"ahler-Weyl connection, then $DK=-\twK\iden+\frac12\dvK J+\frac12(d^DK)^\sd$
where $(d^DK)^\sd$ is a selfdual $2$-form and $\twK$,$\dvK$ are functions.

Now let $\tw=\twK |K|^{-1}$, $\dv=\dvK |K|^{-1}$, $\xi=K|K|^{-1}$,
$\chi=J\xi$. Since $d^DK=\dvK J+(d^DK)^\sd$, it follows that
\begin{align*}
d^DK(K,.)|K|^{-2}&=\dv\chi+(d^DK)^\sd(K,.)|K|^{-2}\\
\tag*{and}(*d^DK)(K,.)|K|^{-2}&=-\dv\chi+(d^DK)^\sd(K,.)|K|^{-2}.
\end{align*}
Therefore $D^{sd}=D+\tw\xi+\dv\chi$ and $(d^DK)^\sd(K,.)|K|^{-2}
=\dv\chi-\omega$.

\begin{thm}\label{csthm} Let $\Mn$ be an oriented conformal $4$-manifold with
conformal vector field $K$ and suppose that $J$ is an invariant antiselfdual
almost complex structure on $\Mn$. Then $J$ is integrable iff
$\chi=J\xi=JK/|K|$ is a shear-free geodesic congruence on the Jones-Tod Weyl
space $B$. Furthermore, the K\"ahler-Weyl structure associated to $J$ is given
by $D=D^{sd}-\tw\xi-\dv\chi$ where
$D^B\chi=\dv(\iden-\chi\tens\chi)+\tw\,{*\chi}$ on $B$.

\proofof{thm} Clearly $\chi$ is invariant and horizontal, hence basic. Let
$\dv,\tw$ be invariant sections of $L^{-1}$ and set
$D=D^{sd}-\tw\xi-\dv\chi$. If $J$ is integrable then we have seen above that
the K\"ahler-Weyl connection is of this form. Therefore it suffices to prove
that $DJ=0$ iff $D^B\chi=\dv(\iden-\chi\tens\chi)+\tw\,{*\chi}$ on $B$.
Since $J=\xi\wedge\chi-*(\xi\wedge\chi)$ this is a straightforward
computation. Let $X$ be any vector field on $\Mn$. Then
\begin{align*}
\D_XJ&=\D_X\xi\wedge\chi+\xi\wedge \D_X\chi
-*(\D_X\xi\wedge\chi-\xi\wedge \D_X\chi).\\
\intertext{Now $D=D^{|K|}+\frac12\omega-\tw\xi-\dv\chi$ and so, since
$D^{|K|}\xi=-\frac12{*\xi\wedge\omega}$, we have}
\D_X\xi&=-\tfrac12\,{*(X\wedge\xi\wedge\omega)}-\tfrac12\cip{\xi,X}\omega
-\tw\bigl(X-\cip{\xi,X}\xi\bigr)+\dv\cip{\xi,X}\chi.\\
\intertext{Also $D=D^B-\frac12\omega-\tw\xi-\dv\chi$ and so}
\D_X\chi&=D^B_X\chi-\tfrac12\omega(\chi)X+\cip{\chi,X}\omega
-\dv\bigl(X-\cip{\chi,X}\chi\bigr)+\tw\cip{\chi,X}\xi.\\
\intertext{Therefore}
\D_X\xi\wedge\chi&=\tfrac12\cip{\chi,X}\,{*(\xi\wedge\omega)}
-\tfrac12\omega(\chi)\,{*(\xi\wedge X)}
-\tw\bigl(X-\cip{\xi,X}\xi\bigr)\wedge\chi\\
&\qquad\qquad-\tfrac12\cip{\xi,X}\omega\wedge\chi\\
\xi\wedge \D_X\chi&=\xi\wedge D^B_X\chi-\tfrac12\omega(\chi)\xi\wedge X
+\tfrac12\cip{\chi,X}\xi\wedge\omega-\dv\xi\wedge\bigl(X-\cip{\chi,X}\chi\bigr)
\end{align*}
and so
\begin{multline*}\xi\wedge \D_X\chi-*(\D_X\xi\wedge\chi)=\\
\xi\wedge D^B_X\chi-\dv\xi\wedge\bigl(X-\cip{\chi,X}\chi\bigr)
+\tw\,{*\bigl((X-\cip{\xi,X}\xi)\wedge\chi\bigr)}
+\tfrac12\cip{\xi,X}\,{*(\omega\wedge\chi)}.
\end{multline*}
Since the right hand side is vertical, it follows that $\D_XJ=0$ iff
$$D^B_X\chi-\cip{D^B_X\chi,\xi}
=\dv\bigl(X-\cip{\chi,X}\chi-\cip{\xi,X}\xi\bigr)
+\tw\,\iota^{\strut}_X\,{*^{\strut}_B\chi}
-\tfrac12\cip{\xi,X}\,{*(\xi\wedge\omega\wedge\chi)}.$$
If $X$ is parallel to $\xi$, this holds automatically
since $\cL_K\chi=0$, and so by considering $X\perp\xi$
we obtain the theorem.
\end{proof}

When $M$ is selfdual, this theorem unifies (the local aspects of) LeBrun's
treatment of scalar-flat K\"ahler metrics with symmetry~\cite{LeBrun1,LeBrun2}
and the hypercomplex structures with symmetry studied by Chave, Tod and
Valent~\cite{CTV} and Gauduchon and Tod~\cite{GT}. To see this, note that
since $D$ is canonically determined by $\Omega_J$, and $\cL_K\Omega_J=0$, it
follows that $\cL_KD=0$ on $L^1$, which means that $d\twK=F^D(K,.)$. Since $K$
is a conformal vector field, it follows that $\cL_KD=0$ on $T\Mn$ as well,
which gives:
$$\tfrac12d\dvK(X) J+\tfrac12\D_X(d^DK)^\sd +W\s/{K,X}-r^D(K)\skwend X
+r^D(X)\skwend K=0.$$
If we contract this with $J$, we obtain $d\dvK=2r^D_0(JK,.)=-i\rho^D(K,.)$.
Thus $\rho^D$ and $F^D$ are the selfdual Maxwell fields associated to the
monopoles $i\dv$ and $\tw$ respectively. Since they are selfdual, it follows
that $d\twK=0$ iff $\Mn,J$ is locally scalar-flat K\"ahler, while $d\dvK=0$ iff
$\Mn,J$ is locally hypercomplex.

Now suppose that $B$ is Einstein-Weyl and that $\gmw$ is any nonvanishing
monopole, and let $\Mn$ be the corresponding selfdual conformal
$4$-manifold. Then \emph{each} shear-free geodesic congruence $\chi$ induces
on $\Mn$ an invariant antiselfdual complex structure $J$.  On the other hand
if we fix $\chi$, then, as we have seen, its divergence and twist, $\dv$ and
$\tw$, are monopoles on $B$.  Using these we can characterise special cases of
the construction as follows.
\begin{enumerate}
\item $(\Mn,J)$ is locally scalar-flat K\"ahler iff
$\tw=a\gmw$ for some constant $a$, and if $a$ is nonzero, we may assume
$a=1$, by normalising $\gmw$.
\begin{itemize}
\item If $\tw=0$ then $\Mn$ is locally scalar-flat K\"ahler
and $K$ is a holomorphic Killing field. If $\dv=b\gmw$, then
$\Mn$ is locally hyperK\"ahler.~\cite{LeBrun1,LeBrun2}
\item If $\tw=\gmw$ then $\Mn$ is locally scalar-flat K\"ahler
and $K$ is a holomorphic homothetic vector field.
\end{itemize}
\item $(\Mn,J)$ is locally hypercomplex iff
$\dv=b\gmw$ for some constant $b$, and if $b$ is nonzero, we may assume
$b=1$, by normalising $\gmw$.
\begin{itemize}
\item If $\dv=0$ then $\Mn$ is locally hypercomplex
and $K$ is a triholomorphic vector field. If $\tw=a\gmw$, then
$\Mn$ is locally hyperK\"ahler.~\cite{CTV,GT}
\item If $\dv=\gmw$ then $\Mn$ is locally hypercomplex
and $K$ is a hypercomplex vector field.
\end{itemize}
\end{enumerate}
Here we say a conformal vector field on a hypercomplex $4$-manifold
is \emphdef{hypercomplex} iff $\cL_KD=0$ where $D$ is the
Obata connection. It follows that for each of the hypercomplex
structures $I$, $\cL_KI$ is a $D$-parallel antiselfdual endomorphism
anticommuting with $I$. The map $I\mapsto\cL_KI\perp I$ is therefore
given by $I\mapsto [cJ,I]$ for one of the hypercomplex structures $J$
and a real constant $c$. Consequently $K$ is holomorphic with respect to
$\pm J$, and is \emphdef{triholomorphic} iff $c=0$.

The twistorial interpretation of the above special cases is as follows.
Firstly, if $\tw=0$ on $B$ then the corresponding line bundle on $\cS$
is trivial; hence so is its pullback to $Z$. On the other hand, if
$\tw=\gmw$ then the line bundle on $\cS$ is nontrivial, but we are
pulling it back to (an open subset of) its total space. Such a pullback
has a tautological section, and hence is trivial away from the
zero section. The story for $\dv$ is similar.

We now combine these observations with the mini-Kerr theorem.

\begin{thm}\label{quotthm}
Let $B$ be an arbitrary three dimensional Einstein-Weyl space.
\begin{enumerate}
\item $B$ may be obtained \textup(locally\textup) as the quotient of a
scalar-flat K\"ahler $4$-manifold by a holomorphic homothetic vector field.
\item It may also be obtained as the quotient of a hypercomplex $4$-manifold
by a hypercomplex vector field.
\item $B$ is locally the quotient of a hyperK\"ahler $4$-manifold by a
holomorphic homothetic vector field if and only if it admits a shear-free
geodesic congruence with linearly dependent divergence and twist.
\end{enumerate}
\proofof{thm} By the mini-Kerr theorem $B$ admits a shear-free geodesic
congruence. The divergence $\dv$ and twist $\tw$ are monopoles on $B$, which
may be used to construct the desired hypercomplex and scalar-flat K\"ahler
spaces wherever they are nonvanishing. The hyperK\"ahler case was
characterised above by the constancy of $\dvK$ and $\twK$. On $B$, this
implies that $\dv$ and $\tw$ are linearly dependent, i.e., $c_1\dv+c_2\tw=0$
for constants $c_1$ and $c_2$. Conversely given an Einstein-Weyl space with a
shear-free geodesic congruence $\chi$ whose divergence and twist satisfy this
condition, any nonvanishing monopole $\gmw$ with $\tw=a\gmw$ and $\dv=b\gmw$
gives rise to a hyperK\"ahler metric (and this $\gmw$ is unique up to a
constant multiple unless $\dv=\tw=0$).
\end{proof}
Maciej Dunajski and Paul Tod~\cite{DT} have recently obtained a related
description of hyperK\"ahler metrics with homothetic vector fields by reducing
Plebanski's equations.

The following diagram conveniently summarises the various Weyl derivatives
involved in the constructions of this section, together with the $1$-forms
translating between them.
\begin{diagram}[width=4.2em,height=2em]
D^{|K|}&\rTo^{+\frac12\omega}&D^{sd}   &\rTo^{+\frac12\omega}&D^B\\
       &\rdTo        &\uTo_{+\tw\xi+\dv\chi}&       &\uTo_{+\tw\xi+\dv\chi}\\
       &                     &D        &\rTo^{+\frac12\omega}&D\S/\chi\\
       &             &                      &\rdTo  &\uTo_{+\tw\xi+\dv\chi}\\
       &                     &         &                     &D^{LW}
\end{diagram}
The Weyl derivatives in the right hand column are so labelled because on $B$
we have $D^{LW}\rTo^{+\dv\chi}D\S/\chi\rTo^{+\dv\chi}D^B$, where $D^B$ is
Einstein-Weyl, $D\S/\chi$ is the Weyl derivative canonically associated to the
congruence $\chi$, and, \emph{in the case that $\tw=0$}, $D^{LW}$ is the
LeBrun-Ward gauge. The central role played by $D^{sd}$ in these constructions
explains the frequent occurrence of the Ansatz $g=Vg\s/B+V^{-1}(dt+A)^2$ for
selfdual metrics with symmetry. In particular, if $g\s/B$ is the LeBrun-Ward
gauge of a LeBrun-Ward geometry and $V$ is a monopole in this gauge, then $g$
is a scalar-flat K\"ahler metric.

\section{Selfdual Einstein $4$-manifolds with symmetry}\label{SDE}

In this section we combine results of Tod~\cite{Tod3} and Pedersen and
Tod~\cite{PT2} to show that the constructions of the previous section
cover essentially all selfdual Einstein metrics with symmetry.
\begin{prop}\label{SymEin}\tcite{PT2} Let $g$ be a four dimensional Einstein
metric with a conformal vector field $K$. Then one of the following must hold:
\begin{enumerate}
\item $K$ is a Killing field of $g$
\item $g$ is Ricci-flat and $K$ is a homothetic vector field \textup(i.e.,
$\cL_K D^g=0$\textup)
\item $g$ is conformally flat.
\end{enumerate}
\end{prop}
Now suppose $g$ is a selfdual Einstein metric with nonzero scalar curvature
and a conformal vector field $K$. Then, except in the conformally flat
case, $K$ is a Killing field of $g$ and so we may apply the following.
\begin{thm}\label{KillEin}\tcite{Tod3} Let $g$ be a selfdual Einstein metric
with nonzero scalar curvature and $K$ a Killing field of $g$. Then the
antiselfdual part of $D^gK$ is nonzero, and is a pointwise multiple of an
integrable complex structure $J$. The corresponding K\"ahler-Weyl structure is
K\"ahler, and $K$ is also a Killing field for the K\"ahler metric.
\end{thm}
If, on the other hand, $\scal^g$ is zero, then $g$ itself is (locally) a
hyperK\"ahler metric and, unless $g$ is conformally flat, $\cL_KD^g=0$, and so
$K$ is a hypercomplex vector field. In the conformally flat case, $K$ may not
be a homothety of $g$, but it is at least a homothety with respect to
\emph{some} compatible flat metric. Thus, in any case, the conformal vector
field $K$ is holomorphic with respect to some K\"ahler structure on $\Mn$.

We end this section by noting that in the case of selfdual Einstein metrics
with Killing fields, Tod's work~\cite{Tod3} shows how to recover the Einstein
metric from the LeBrun-Ward geometry. More precisely, if $M$ is a selfdual
Einstein $4$-manifold with a Killing field, and $B$ is the LeBrun-Ward
quotient of the corresponding scalar-flat K\"ahler metric, then
either $B$ is flat, or the monopole defining $M$ is of the form
$$\gmw=\bigl(a(1-\tfrac12zu_z)+\tfrac12 bu_z\bigr)\mu_{LW}^{-1},$$ where
$u(x,y,z)$ is the solution of the $\SU(\infty)$ Toda field equation, and
$a,b\in\R$ are not both zero. Conversely, for any LeBrun-Ward geometry (given
by $u$), the section $\bigl(a(1-\frac12zu_z)+\frac12bu_z\bigr)\mu_{LW}^{-1}$
of $L^{-1}$ is a monopole for any $a,b\in\R$, and if $g\s/K$ is the
corresponding K\"ahler metric, then $(az-b)^{-2}g\s/K$ is Einstein with scalar
curvature $-12a$. When $a=0$, we reobtain the case of hyperK\"ahler metrics
with Killing fields, while if $a\neq0$, one can set $b=0$ by translating the
$z$ coordinate (although $u$ will be a different function of the new $z$
coordinate).

\section{Einstein-Weyl structures from $\R^4$}\label{EWR4}

Our aim in the remaining sections is to unify and extend many of the examples
of K\"ahler-Weyl structures with symmetry studied up to the present, using the
framework developed in sections~\ref{CSKW}--\ref{JTsect}. We discuss
both the simplest and most well known cases and also more complicated examples
which we believe are new. We begin with $\R^4$.

A conformally flat $4$-manifold is both selfdual and antiselfdual, so when we
apply the Jones and Tod construction we have the freedom to reverse the
orientation. Consequently, not only is $D^B=D^{|K|}+\omega$ Einstein-Weyl, but
so is $\tilde D^B=D^{|K|}-\omega$. Therefore
$0=\sym_0(D^B\omega+\omega\tens\omega)=\sym_0 D^{|K|}\omega =\sym_0(\tilde
D^B\omega-\omega\tens\omega)$.  Since $|K|^{-1}$ is a monopole,
$g=|K|^{-2}\conf\s/B$ (the gauge in which the monopole is constant) is a
\emphdef{Gauduchon metric} in the sense that $\omega$ is divergence-free with
respect to $D^g=D^{|K|}$. It follows that $\omega$ is dual to a Killing
field of $g$. Furthermore, the converse is also true: that is, if
$D^B=D^g+\omega$ is Einstein-Weyl and $\omega$ is dual to a Killing field of
$g$, then $\tilde D^B=D^g-\omega$ is also Einstein-Weyl, and therefore the
$4$-manifold $\Mn$ given by the monopole $\mu_g^{-1}$ is both selfdual and
antiselfdual, hence conformally flat.

The condition that an Einstein-Weyl space admits a compatible metric $g$ such
that $D=D^g+\omega$ with $\omega$ dual to a Killing field of $g$ is of
particular importance because it always holds in the compact case: on any
compact Weyl space there is a Gauduchon metric $g$ unique up to
homothety~\cite{Gauduchon1}, and $g$ has this additional property when the
Weyl structure is Einstein-Weyl~\cite{Tod1}. Consequently, the local
quotients of conformally flat $4$-manifolds exhaust the possible local
geometries of compact Einstein-Weyl $3$-manifolds. These geometries were
obtained in~\cite{PT1} as local quotients of $S^4$. Now any conformal vector
field $K$ on $S^4$ has a zero and is a homothetic vector field with respect to
the flat metric on $\R^4$ given by stereographic projection away from any such
zero. Hence we can view these Einstein-Weyl geometries as local quotients of
the flat metric on $\R^4$ by a homothetic vector field and use the
constructions of section~\ref{JTsect} to understand some of their properties.

Suppose first that $K$ vanishes on $R^4$ and let the origin be such a zero.
Then $K$ generates one parameter group of \emphdef{linear} conformal
transformations of the flat metric $g$. This is case 1 of~\cite{PT1} and
we may choose coordinates such that
\begin{align*}
g&=dr^2+\tfrac14r^2\bigl(d\theta^2+\sin^2\theta\,d\phi^2
+(d\psi+\cos\theta\,d\phi)^2\bigr)\\
K&=a r\frac\partial{\partial r}-(b+c)\frac\partial{\partial\phi}
-(b-c)\frac\partial{\partial\psi}.
\end{align*}
Note that $K$ is also a homothety of the flat metric $\tilde g=r^{-4}g$
obtained from $g$ by the orientation reversing conformal transformation
$r\mapsto\tilde r=1/r$. With a fixed orientation,
\begin{align*}
D^gK&=\hphantom{-}a\,\iden+\tfrac12(b+c)J^\sd+\tfrac12(b-c)J^\asd\\
D^{\tilde g}K&=-a\,\iden+\tfrac12(b-c)\tilde J^\sd+\tfrac12(b+c)\tilde J^\asd
\end{align*}
where $J^\sdasd$ are $D^g$-parallel complex structures on $\R^4$, one
selfdual, the other antiselfdual, and, similarly, $\tilde J^\sdasd$ are
$D^{\tilde g}$-parallel. The Weyl structures $D^{|K|}\pm\omega$ are
Einstein-Weyl on the quotient $B$, where
$$\omega=\frac{(b+c)g(J^\sd K,.)-(b-c)g(J^\asd K,.)}{g(K,K)}
=\frac{(b-c)\tilde g(\tilde J^\sd K,.)-(b+c)\tilde g(\tilde J^\asd K,.)}
{\tilde g(K,K)}.$$
Without loss of generality, we consider only $D^B=D^{|K|}+\omega$.
By Theorem~\ref{csthm}, $J^\asd K$ and $\tilde J^\asd K$
generate shear-free geodesic congruences with $\dv^\asd=(b+c)|K|^{-1}$,
$\tilde\dv^\asd=(b-c)|K|^{-1}$ and $\tw^\asd=a|K|^{-1}=-\tilde\tw^\asd$.

If $b^2=c^2$, then $K$ is triholomorphic, and so the quotient geometry is
hyperCR: it is the Berger sphere family. If we take $b=c$ then $J^\asd$ is no
longer unique, and the hyperCR structure is given by the congruences
associated to $JK$, where $J$ ranges over the parallel antiselfdual
complex structures of $g$; $\tilde J^\asd K$, by contrast, is the geodesic
symmetry $\partial/\partial\phi$ of $B$. In addition, the antiselfdual
rotations all commute with $K$, so $B$ has a four dimensional symmetry group,
locally isomorphic to $S^1\times S^3$.

If $bc=0$, then although $K$ is not a Killing field on $\R^4$ unless $a=0$, it
is Killing with respect to the product metric on $S^2\times\cH^2$ which is
scalar flat K\"ahler (where the hyperbolic metric on $\cH^2$ has equal and
opposite curvature to the round metric on $S^2$) and conformal to
$\R^4\setdif\R$. Hence these quotients are Toda.

If $a=0$, then $K$ is a Killing field, and so the (local) quotient geometry
is also Toda, simply because it is the quotient of a flat metric by a Killing
field.

If $b^2=c^2$ and $bc=0$ then $b=c=0$ and the quotient is the round $3$-sphere,
while if $a=0$ and $bc=0$ it turns out to be the hyperbolic metric. If $a=0$
and $b^2=c^2$, the quotient geometry is the flat Weyl space: the hyperCR
congruences become the translational symmetries, and (for $b=c$) $\tilde
J^\asd K$ is the radial symmetry.

We now briefly consider the case that $K$ does not vanish on $\R^4$ (and so is
not linear with respect to any choice of origin). This is case 2
of~\cite{PT1}, and we may choose a flat metric $g$ with respect to which $K$
is a transrotation. Since $K$ is a Killing field, the quotient Einstein-Weyl
space is Toda. For $b=0$, it is flat, while for $c=0$ we obtain $\cH^3$.

\section{HyperK\"ahler metrics with triholomorphic Killing fields}\label{HKt}

If $M$ is a hyperK\"ahler $4$-manifold and $K$ is a triholomorphic Killing
field, then $\dv$ and $\tw$ both vanish, so the corresponding Einstein-Weyl
space is flat and the congruence consists of parallel straight lines.
HyperK\"ahler $4$-manifolds with triholomorphic Killing fields therefore
correspond to nonvanishing solutions of the Laplace equation on an open subset
of $\R^3$, or some discrete quotient. This is the Gibbons-Hawking Ansatz for
selfdual Euclidean vacua~\cite{GH}.

In~\cite{Ward}, Ward used this Ansatz to generate new Toda Einstein-Weyl
spaces from axially symmetric harmonic functions. The idea is beautifully
simple: since the harmonic function is preserved by a Killing field on $\R^3$,
the Gibbons-Hawking metric admits a two dimensional family of commuting
Killing fields; one of these is triholomorphic, but the others need not be,
and so they have other Toda Einstein-Weyl spaces with symmetry as quotients.

Let us carry out this procedure explicitly. In cylindrical polar coordinates
$(\eta,\rho,\phi)$, the flat metric is $d\eta^2+d\rho^2+\rho^2d\phi^2$ and
the generator of the axial symmetry is $\partial/\partial\phi$. An invariant
monopole (in the gauge determined by the flat metric) is a function
$W(\rho,\eta)$ satisfying $\rho^{-1}(\rho W_\rho)_\rho+W_{\eta\eta}=0$.  Note
that if $W$ is a solution of this equation, then so is $W_\eta$, and $W_\eta$
determines $W$ up to the addition of $C_1\log(C_2\rho)$ for some
$C_1,C_2\in\R$. This provides a way of integrating the equation $d{*dW}=0$ to
give $*dW=dA$: if we write $W=V_\eta$, then we can take $A=\rho
V_\rho\,d\phi$. This choice of integral determines the lift of
$\partial/\partial\phi$ to the $4$-manifold. The hyperK\"ahler metric is
$$g=V_\eta(d\eta^2+d\rho^2+\rho^2d\phi^2)
+V_\eta^{-1}(d\psi+\rho V_\rho\,d\phi)^2.$$
In order to take the quotient by $\partial/\partial\phi$, we rediagonalise:
$$g=V_\eta\biggl(d\rho^2+d\eta^2 +\frac1{V_\eta^2+V_\rho^2}d\psi^2\biggr)
+\frac{\rho^2(V_\eta^2+V_\rho^2)}{V_\eta}
\biggl(d\phi+\frac{V_\rho}{\rho(V_\eta^2+V_\rho^2)}d\psi\biggr)^2.$$
We now recall that the hyperK\"ahler metric lies midway between the constant
length gauge of $\partial/\partial\phi$ and the LeBrun-Ward gauge of the
quotient. Consequently we find that $D^B=D^{LW}+\omega$ where:
\begin{align*}
g\s/{LW}&=\rho^2(V_\eta^2+V_\rho^2)(d\rho^2+d\eta^2)+\rho^2d\psi^2
=\rho^2(dV^2+d\psi^2)+(\rho V_\rho\, d\eta-\rho V_\eta\,d\rho)^2\\
\omega&=\frac{2V_\eta}{\rho^2(V_\eta^2+V_\rho^2)}
(\rho V_\rho\, d\eta-\rho V_\eta\,d\rho).
\end{align*}
Note that $d(\rho V_\rho d\eta-\rho V_\eta d\rho)=0$. This can be
integrated by writing $V=U_\eta$, with $U(\rho,\eta)$ harmonic. Then $z=\rho
U_\rho$ parameterises the hypersurfaces orthogonal to the shear-free
twist-free congruence, and isothermal coordinates on these hypersurfaces are
given by $x=U_\eta$, $y=\psi$. Hence, although the Einstein-Weyl space is
completely explicit, the solution $e^u=\rho^2$ of the $\SU(\infty)$ Toda field
equation is only given implicitly. Nevertheless, we \emph{have} found the
congruence, the isothermal coordinates and the monopole $u_z\mu_{LW}^{-1}$.

The symmetry $\partial/\partial\psi$, like the axial symmetry
$\partial/\partial\phi$ on $\R^3$, generates a congruence which is
divergence-free and twist-free, although it is not geodesic. For this reason
it is natural to say that Ward's spaces are Einstein-Weyl \emphdef{with an
axial symmetry}. They are studied in more detail in~\cite{DMJC4}.

\medbreak

Ward's construction can be considerably generalised. First of all, one can
obtain new Toda Einstein-Weyl spaces by considering harmonic functions
invariant under other Killing fields. The general Killing field on $\R^3$ may
be taken, in suitably chosen cylindrical coordinates, to be of the form
$b\partial/\partial\phi+c\partial/\partial\eta$ for $b,c\in\R$. By introducing
new coordinates $\zeta=(b\eta-c\phi)/\sqrt{b^2+c^2}$ and
$\theta=(b\phi+c\eta)/\sqrt{b^2+c^2}$, so that the Killing field is a multiple
of $\partial/\partial\theta$, one can carry out the same procedure as before
to obtain the following Toda Einstein-Weyl spaces:
\begin{align*}
g\s/{LW}&=G(\rho,\zeta)\bigl(d\rho^2+F(\rho)d\zeta^2\bigr)
+\rho^2F(\rho)^{-1}\beta^2\\
&=\rho^2\biggl(dV^2+\frac1{b^2+c^2}\Bigl(
c\bigl[\rho V_\rho\, d\zeta-F(\rho)^{-1}\rho V_\zeta\,d\rho\bigr]
+b\,d\psi\Bigr)^2\biggr)\\
&\qquad\qquad\qquad\qquad+\frac1{b^2+c^2}
\Bigl(b\bigl[\rho V_\rho\, d\zeta-F(\rho)^{-1}\rho V_\zeta\,d\rho\bigr]
-c\,d\psi\Bigr)^2\\
\omega&=\frac{2bV_\zeta}{(b^2+c^2)G(\rho,\zeta)}
\Bigl(b\bigl[\rho V_\rho\, d\zeta-F(\rho)^{-1}\rho V_\zeta\,d\rho\bigr]
-c\,d\psi\Bigr),
\intertext{where}
F(\rho)&=\frac{(b^2+c^2)\rho^2}{b^2\rho^2+c^2},\qquad
G(\rho,\zeta)=\frac{(b^2\rho^2+c^2)V_\zeta^2+(b^2+c^2)\rho^2V_\rho^2}
{b^2+c^2},\\
\tag*{and}\beta&=
d\psi-\frac{bc(1-\rho^2)}{b^2\rho^2+c^2}\bigl[\rho V_\rho\,d\eta
-F(\rho)^{-1}\rho V_\eta\,d\rho\bigr].
\end{align*}
Note that the symmetry $\partial/\partial\psi$ is twist-free if and only
if $bc=0$. When $b=0$, the Toda Einstein-Weyl space is just $\R^3$ (the
only Einstein-Weyl space with a parallel symmetry),
while $c=0$ is Ward's case.

\medbreak

A further generalisation of this procedure is obtained by observing that the
flat Weyl structure on $\R^3$ is preserved not just by Killing fields, but by
homothetic vector fields. Now, for a section $\gmw$ of $L^{-1}$, invariance no
longer means that the function $\gmw\mu_{\R^3}$ is constant along the flow of
the homothetic vector field, since the length scale $\mu_{\R^3}$ is not
invariant. Hence it is better to work in a gauge in which the homothetic
vector field is Killing. To do this we may choose spherical polar coordinates
$(r,\theta,\phi)$ such that the flat Weyl structure on $\R^3$ is
\begin{align*}
g_0&=r^{-2}dr^2+d\theta^2+\sin^2\theta\,d\phi^2\\
\omega_0&=r^{-1}dr
\end{align*}
and the homothetic vector field is a linear combination of $r\partial/\partial
r$ and $\partial/\partial\phi$. For simplicity, we shall only consider here
the case of a pure dilation $X=r\partial/\partial r$. If $\gmw=W\mu_0^{-1}$ is
an invariant monopole (where $\mu_0$ is the length scale of $g_0$) then
$W_r=0$ and $W(\theta,\phi)$ is a harmonic function on $S^2$. We write
$g\s/{S^2}=\sigma_1^2+\sigma_2^2$ and $W=\frac12(h+\overline h)$ with $h$
holomorphic on an open subset of $S^2$. Then the hyperK\"ahler metric is
$$g=\frac{r(h+\overline h)}{2|h|^2}
\bigl(|h|^2(\sigma_1^2+\sigma_2^2)+\beta^2\bigr)
+\frac{2|h|^2}{(h+\overline h)r}\bigl(dr+i(h+\overline h)r\,\beta\bigr)^2,$$
where $\beta$ is a $1$-form on $S^2$ with $d\beta=\frac12(h+\overline h)
\sigma_1\wedge\sigma_2$. One easily verifies that the quotient space is the
Einstein-Weyl space with geodesic symmetry given by the holomorphic function
$H=1/h$.

The computation for the general homothetic vector field is more complicated,
but one obtains Gibbons-Hawking metrics admitting holomorphic conformal vector
fields which are neither triholomorphic or Killing, and therefore, as
quotients, explicit examples of Einstein-Weyl spaces (with symmetry) which are
neither hyperCR, nor Toda, yet they admit a shear-free geodesic congruence
with linearly dependent divergence and twist.

\section{Congruences and monopoles on $\cH^3$, $\R^3$ and $S^3$}\label{CMCC}

An important special case of the theory presented in this paper is the case
of monopoles on spaces of constant curvature. Since each shear-free geodesic
congruence on these spaces induces a complex structure on the selfdual space
associated to any monopole, it is interesting to find such congruences.

The twist-free case has been considered by Tod~\cite{Tod2}. In this case we
have a LeBrun-Ward space of constant curvature, given by a solution $u$ of
the Toda field equation with $u_z\,dz$ exact. This happens precisely when
$u(x,y,z)=v(x,y)+w(z)$. The solutions, up to changes of isothermal
coordinates, are given by
$$e^u=\frac{4(az^2+bz+c)}{(1+a(x^2+y^2))^2}$$
where $a,b,c$ are constants constrained by positivity. As shown
in~\cite{Tod2}, there are essentially six cases: three on hyperbolic space
($b^2-4ac>0$), two in flat space ($b^2-4ac=0$), and one on the sphere
($b^2-4ac<0$). One of the congruences in each case
is a radial congruence, orthogonal to distance spheres. The other two types of
congruences on hyperbolic space are orthogonal to horospheres and hyperbolic
discs respectively, while the other type of congruence on flat space is
translational. Only the radial congruences have singularities, and in the flat
case, even the radial congruence is globally defined on $S^1\times S^2$.  We
illustrate the congruences in the following diagrams. 

\bigbreak \ifaddpics
\hspace{3.25cm}$b^2>4ac$\hspace{1.5cm}$b^2=4ac$\hspace{1.5cm}$b^2<4ac$

\bigbreak
$a>0$\hspace{2cm}\getpic{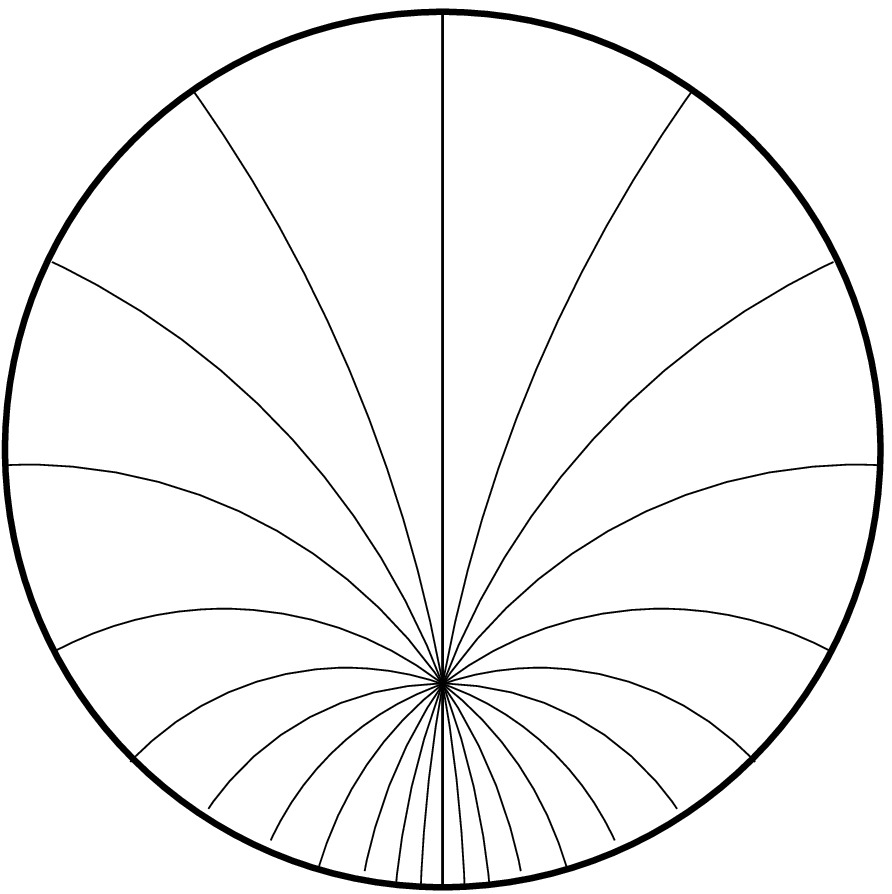}\hspace{1cm}\getpic{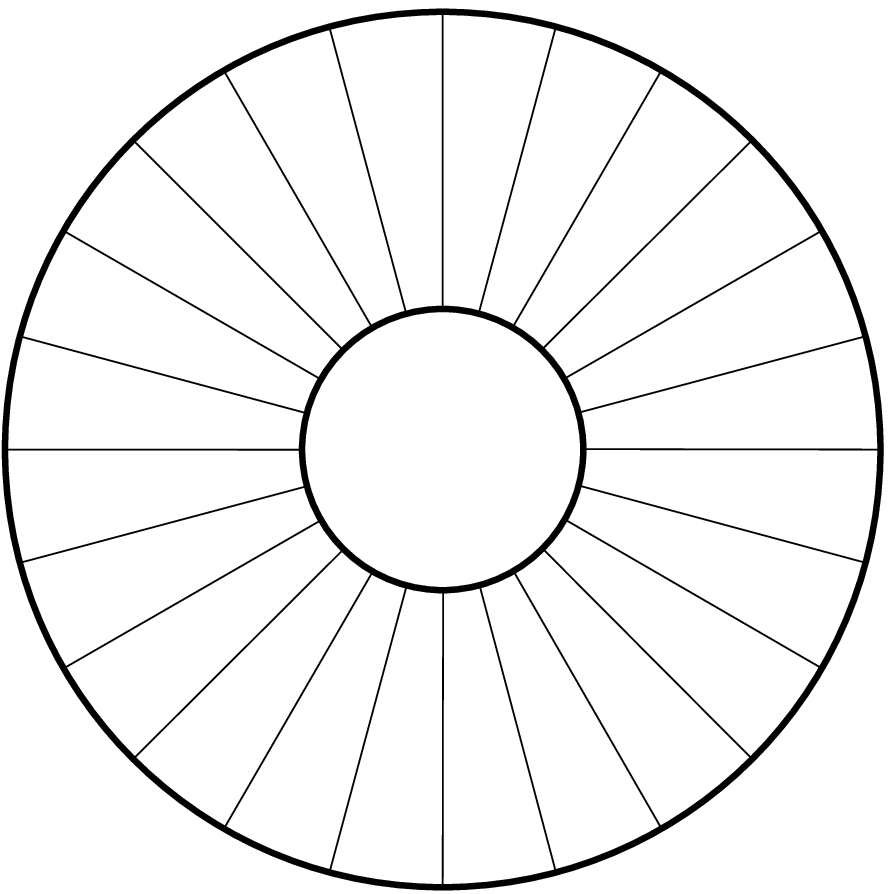}
\hspace{1cm}\getpic{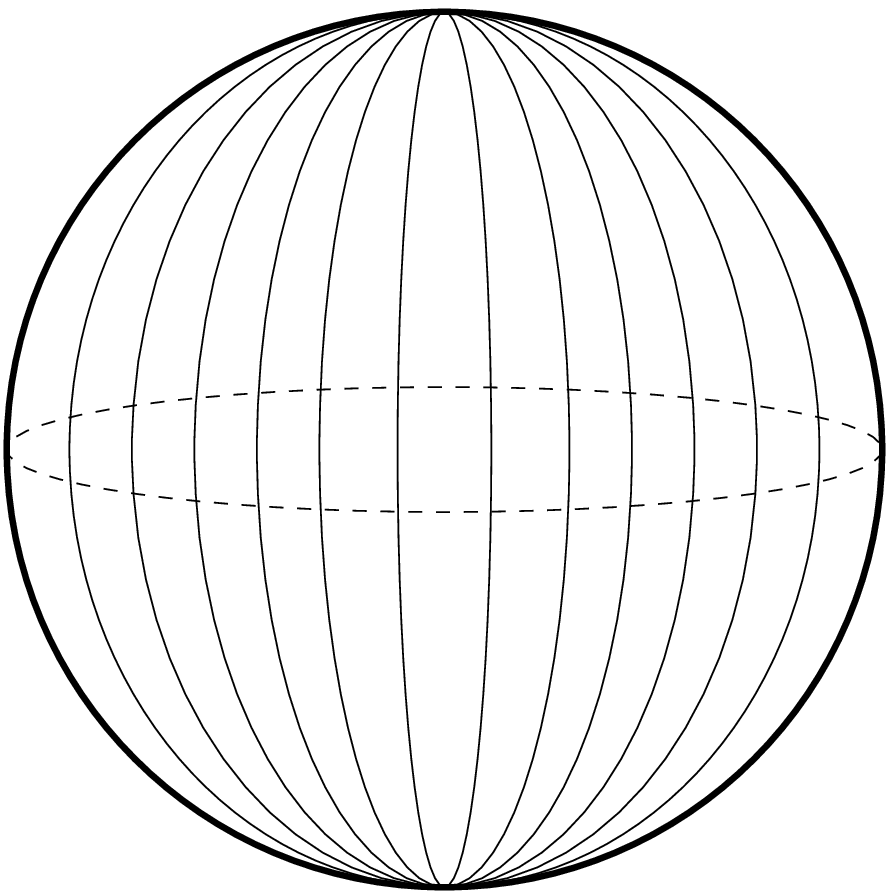}

\medbreak
$a=0$\hspace{2cm}\getpic{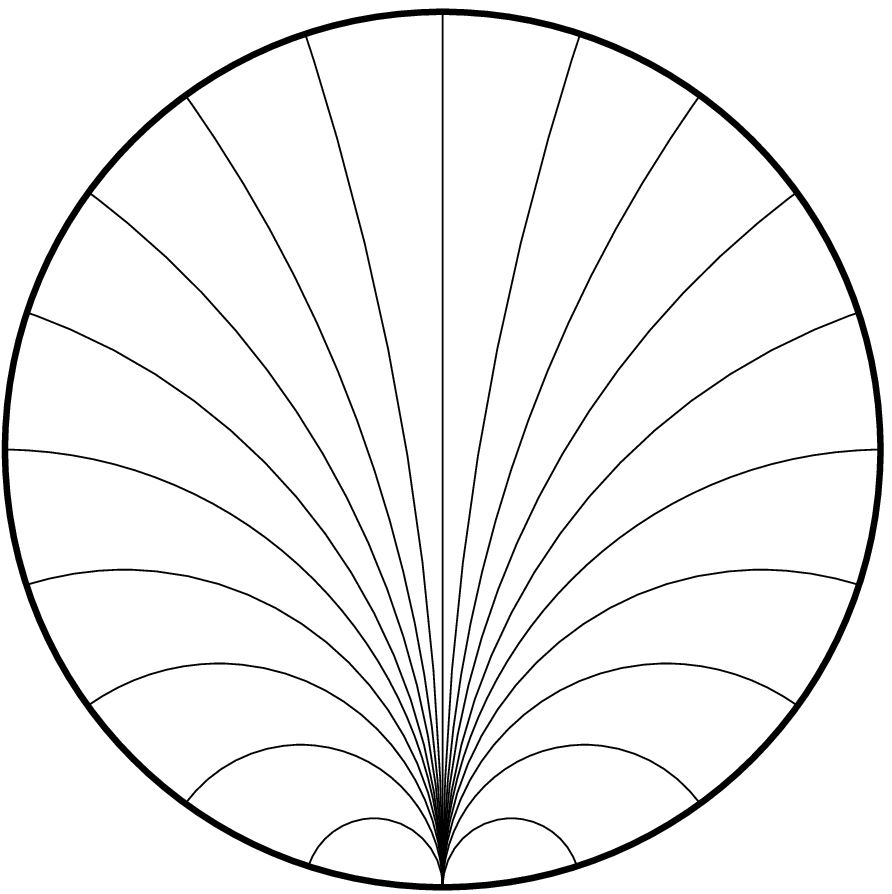}\hspace{1cm}\getpic{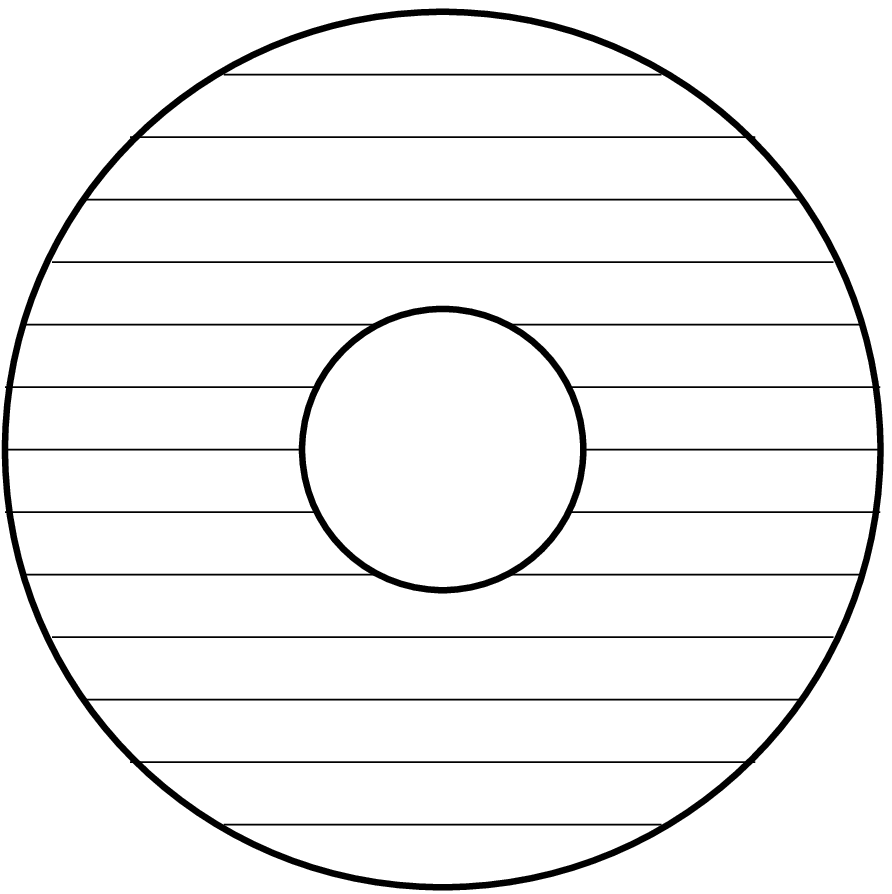}

\medbreak
$a<0$\hspace{2cm}\getpic{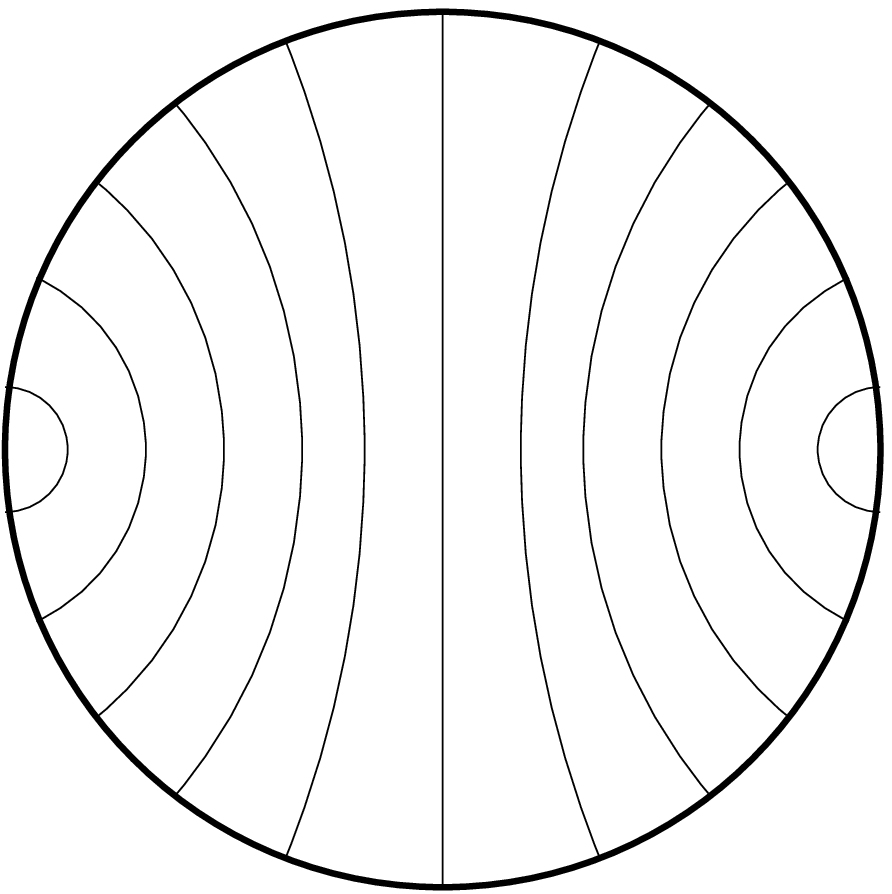}

\bigbreak\fi

The congruences on hyperbolic space $\cH^3$ have been used by LeBrun
(see~\cite{LeBrun1,LeBrun2}) to construct selfdual conformal structures on
complex surfaces. The first type of congruence gives scalar-flat K\"ahler
metrics on blow ups of line bundles over $\CP1$. The second type gives
asymptotically Euclidean scalar-flat K\"ahler metrics on blow-ups of $\C^2$
and hence selfdual conformal structures on $k\CP2$ and closed K\"ahler-Weyl
structures on blow-ups of Hopf surfaces. The final type of congruence descends
to quotients by discrete subgroups of $\SL(2,\R)$ and leads to scalar-flat
K\"ahler metrics on ruled surfaces of genus$\geq2$.

If we look instead for hyperCR structures (i.e., divergence-free congruences),
we have, in addition to the translational congruences on $\R^3$, two such
structures on $S^3$: the left and right invariant congruences, but this
exhausts the examples on spaces of constant curvature. Of course
there is still an abundance of congruences which are neither twist-free nor
divergence-free. For instance, on $\R^3$, a piece of a minitwistor line and
its conjugate define a congruence on some open subset: if the line is real
then this is a radial congruence, but in general, we get a congruence of
rulings of a family of hyperboloids.

\ifaddpics\centerline{\getbigpic{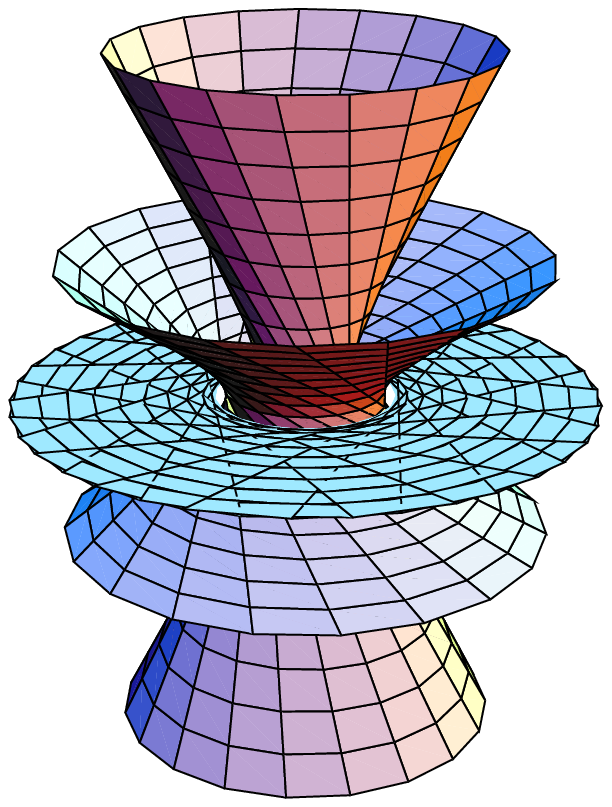}}\fi
This congruence is globally defined on the nontrivial double cover of
$\R^3\setdif S^1$. Its divergence and twist are closely related to the
Eguchi-Hanson I metric as we shall see below.

In general, a holomorphic curve in the minitwistor space of $\R^3$ corresponds
to a null curve in $\C^3$ and the associated congruence consists of the real
points in the tangent lines to the null curve. Since null curves may be
constructed from their real and imaginary parts, which are conjugate minimal
surfaces in $\R^3$, this shows that more complicated congruences are
associated with minimal surfaces.

\medbreak

Turning now to monopoles, we have two simple and explicit types of solutions
of the monopole equation: the constant solutions and the fundamental
solutions. Linear combinations of these give rise to an interesting family of
selfdual conformal structures whose properties are given by the above
congruences. Since such monopoles are spherically symmetric, these selfdual
conformal structures will admit local $\Un(2)$ or $S^1\times\SO(3)$ symmetry.

The $3$-metric with constant curvature $c$ is
$$g_c=\frac4{(1+cr^2)^2}(dr^2+r^2 g\s/{S^2})$$ and the monopoles of interest
are $a+bz$, where $z=(1-cr^2)/2r$ is the fundamental solution centred at
$r=0$. The fundamental solution is the divergence of the radial congruence,
and if we use the coordinate $z$ in place of $r$, we obtain
$$g_c=\biggl(\frac{dz}{z^2+c}\biggr)^2+\frac1{z^2+c}g\s/{S^2}.$$
Rescaling by $(z^2+c)^2$ gives the Toda solution
\begin{align*}
g\s/{LW}&=(z^2+c)g\s/{S^2}+dz^2\\
\omega_{LW}&=-\frac{2z}{z^2+c}dz.
\end{align*}
In the LeBrun-Ward gauge, the monopoles of interest are $\gmw=(a+bz)/(z^2+c)$.
If $c\neq0$ then $\gmw=\frac ac(1-\frac12zu_z)+\frac b2u_z$ and so we may apply
Tod's prescription for the construction of Einstein metrics with symmetry.
Rescaling by $(a^2+c^2)/c$ gives the Einstein metric
$$g=\frac{a^2+c^2}{(az-bc)^2}\biggl(\frac{a+bz}{z^2+c}
\bigl(dz^2+(z^2+c)g\s/{S^2}\bigr) +\frac{z^2+c}{a+bz}(dt+A)^2\biggr)$$ of
scalar curvature $-12ac/(a^2+c^2)$, where $dA=*D\gmw=b\,\mathit{vol}\s/{S^2}$.
This is easily integrated by $A=-b\cos\theta\,d\phi$ where
$g\s/{S^2}=d\theta^2+\sin^2\theta\,d\phi^2$. These metrics are also
well-defined when $c=0$ when they become Taub-NUT metrics with triholomorphic
Killing field $\partial/\partial\psi$. They are also Gibbons-Hawking metrics
for $a=0$, when we obtain the Eguchi-Hanson I and II metrics: this time
$\partial/\partial\phi$ (and the other infinitesimal rotations of $S^2$) is a
triholomorphic Killing field. To relate the metrics to those of~\cite{HP}, one
can set $z=1/\rho^2$ and rescale by a further factor $1/4$.  Then
$$g=\frac{a^2+c^2}{(a-bc\rho^2)^2}\biggl(\frac{a\rho^2+b}{1+c\rho^4}d\rho^2
+\frac14\rho^2\biggl[(a\rho^2+b)g\s/{S^2}
+\frac{1+c\rho^4}{a\rho^2+b}(dt-b\cos\theta\,d\phi)^2\biggr]\biggr)$$
is Einstein with scalar curvature $-48ac/(a^2+c^2)$.
Up to homothety, this is really only a one parameter family of
Einstein metrics, since the original constant curvature metric and
the monopole $\gmw$ can be rescaled. However, the use of three parameters
enables all the limiting cases to be easily found.

These metrics are all conformally scalar-flat K\"ahler via the radial Toda
congruences~\cite{LeBrun0}. The metrics over $\cH^3$ are also conformal to
other scalar-flat K\"ahler metrics, via the horospherical and disc-orthogonal
congruences. The translational congruences on $\R^3$ correspond to the
hyperK\"ahler structures associated with the Ricci-flat $c=0$ metrics. The
metrics coming from $S^3$ admit two hypercomplex structures (coming from the
hyperCR structures), which explains an observation of Madsen~\cite{Madsen}. In
particular when $a=0$, the Eguchi-Hanson I metric has two additional
hypercomplex structures with respect to which $\partial/\partial\psi$ is
triholomorphic. On the other hand, although $\partial/\partial\phi$ is
triholomorphic with respect to the hyperK\"ahler metric, it only preserves one
complex structure from each of these additional families. The corresponding
congruences on $\R^3$ are the two rulings of the families of hyperboloids,
which have the same divergence but opposite twist.  The monopole giving
Eguchi-Hanson I must be the divergence of this congruence.

In~\cite{PT2}, it is claimed that the above constructions give all the
Einstein metrics over $\cH^3$. This is not quite true, because we have not yet
considered the Einstein metrics associated to the horospherical and
disc-orthogonal congruences. These turn out to give Bianchi type VII$_0$ and
VIII analogues of the above Bianchi type IX metrics (by which we mean, the
$\SU(2)$ symmetry group is replaced by Isom($\R^2$) and $\SL(2,\R)$
respectively---see~\cite{Tod2a}). This omission from~\cite{PT2} was simply due
to the nowhere vanishing conformal vector fields on hyperbolic space being
overlooked.

\section{K\"ahler-Weyl spaces with torus symmetry}

On an Einstein-Weyl space with symmetry, an invariant shear-free geodesic
congruence and an invariant monopole together give rise to a selfdual
K\"ahler-Weyl structures, possessing, in general, only two continuous
symmetries. Many explicit examples of such Einstein-Weyl spaces with symmetry
were given in section~\ref{HKt}. Being quotients of Gibbons-Hawking metrics,
these spaces already come with invariant congruences, and solutions of the
monopole equation can be obtained by introducing an additional invariant
harmonic function on $\R^3$, lifting it to the Gibbons-Hawking space, and
pushing it down to the Einstein-Weyl space. Carrying out this procedure in
full generality would take us too far afield, so we confine ourselves to the
two simplest classes of examples: the Einstein-Weyl spaces with axial
symmetry, and the Einstein-Weyl spaces with geodesic symmetry.

We first consider the case of axial symmetry, when the K\"ahler-Weyl structure
is (locally) scalar flat K\"ahler. In~\cite{Joyce}, Joyce constructs such
torus symmetric scalar-flat K\"ahler metrics from a linear equation on
hyperbolic $2$-space. In this way he obtains selfdual conformal structures
on $k\CP2$, generalising (for $k\geq4$) those of LeBrun~\cite{LeBrun1}.
Joyce does not consider the intermediate Einstein-Weyl spaces in his
construction, but one easily sees that his linear equation is equivalent
to the equation for axially symmetric harmonic functions, and that the
associated Einstein-Weyl spaces are precisely the ones with axial
symmetry~\cite{DMJC4}.

Let us turn now to the spaces with geodesic symmetry, where a monopole
invariant under the symmetry is given by a nonvanishing holomorphic function
on an open subset of $S^2$. Indeed, if we write (as before)
\begin{align*}
g&=|H|^{-2}(\sigma_1^2+\sigma_2^2)+\beta^2\\
\omega&=\tfrac i2(H-\overline H)\beta
\end{align*}
with $\beta$ dual to the symmetry, then an invariant monopole in this gauge is
given by the pullback $V$ of a harmonic function on an open subset of $S^2$,
as one readily verifies by direct computation. Hence
$V=\frac12(F+\overline F)$ for some holomorphic function $F$.
The selfdual space constructed from $V$ will admit a K\"ahler-Weyl structure
(coming from the geodesic symmetry) and also a hypercomplex structure (coming
from the hyperCR structure). By Proposition~\ref{gschar}, the geodesic
symmetry preserves the hyperCR congruences, and so it lifts to a
triholomorphic vector field of the hypercomplex structure. Since~\ref{gschar}
is a characterisation, we immediately deduce the following result.

\begin{thm}\label{hcthm} Let $M$ be a hypercomplex $4$-manifold with a two
dimensional family of commuting triholomorphic vector fields. Then the
quotient of $M$ by any of these vector fields is Einstein-Weyl with a geodesic
symmetry, and so the conformal structure on $M$ depends explicitly on two
holomorphic functions of one variable.
\end{thm}

There are two special choices of monopole on such an Einstein-Weyl space: the
$\tw$ and $\dv$ monopoles of the geodesic symmetry. The $\tw$ monopole $(F=H)$
leads us back to the Gibbons-Hawking hyperK\"ahler metric, but the $\dv$
monopole $(F=iH)$ is more interesting. In this case, the K\"ahler-Weyl
structure given by the geodesic symmetry is hypercomplex and so these torus
symmetric selfdual spaces are hypercomplex in two ways. The symmetries are
both triholomorphic with respect to the first hypercomplex structure, but only
one of them is triholomorphic with respect to the additional hypercomplex
structure. If we take the quotient by the bi-triholomorphic symmetry, we
obtain an Einstein-Weyl space with two hyperCR structures, which must be
$S^3$. Hence the spaces with geodesic symmetry, as well as coming from
invariant monopoles on $\R^3$, also come from invariant monopoles on $S^3$.

We end by discussing a third situation in which the spaces with geodesic
symmetry occur. This involves some explicit new solutions~\cite{CT} of the
$\SU(\infty)$ Toda field equation generalising the solutions on $S^3$
described earlier. The corresponding LeBrun-Ward geometries are:
\begin{align*}
g&=(z+h)(z+\overline h)(\sigma_1^2+\sigma_2^2)+dz^2,\\
\omega&=-\frac{2z+h+\overline h}{(z+h)(z+\overline h)}dz,
\end{align*}
where $h$ is an arbitrary nonvanishing holomorphic function on an open subset
of $S^2$. These spaces have no symmetries and so one obtains from them
Einstein metrics with a one dimensional isometry group.  However,
$\partial/\partial z$ does lift to a shear-free congruence on the Einstein
space, and a generalised Jones and Tod construction may be used to show that
the quotient by this conformal submersion is the Einstein-Weyl space with
geodesic symmetry given by $H=1/h$~\cite{CT}. In fact this was how these
interesting Einstein-Weyl spaces were found.

\end{document}